# A Heuristic Framework of Variable Neighborhood Descent Methods for the Large-Scale Multi-Level Facility Location Problem in Supply Chain Networks

Haibo Wang

A. R. Sanchez, Jr. School of Business, Texas A&M International University, Laredo, Texas, USA
hwang@tamiu.edu

Bahram Alidaee

School of Business Administration, The University of Mississippi, University, Mississippi, USA
alidaee@olemiss.edu

# A Heuristic Framework of Variable Neighborhood Descent Methods for the Large-Scale Multi-Level Facility Location Problem in Supply Chain Networks

**Abstract**

This paper addresses the single-assignment, uncapacitated, multi-level facility location (MFL) problem, a strategic decision-making process critical to the design of long-term supply chain networks. Specifically, we examine four- and five-level facility location structures (k-LFL), modeled as a location-allocation problem where demand nodes must be assigned to open facilities across hierarchical levels. Although the MFL has been addressed in the literature, solutions to large-scale, realistic problems involving thousands of nodes are lacking. This paper proposes a heuristic framework based on the Variable Neighborhood Descent (VND) metaheuristic with a multi-start strategy. We develop and compare four variants: Basic Variable Neighborhood Descent (BVND), Pipe Variable Neighborhood Descent (PVND), Cyclic Variable Neighborhood Descent (CVND), and Union Variable Neighborhood Descent (UVND). In each case, a multi-start strategy with strong diversification components is employed. Extensive computational experiments compare the methods on large-scale instances involving up to 10,000 customers, 150 distribution centers, 50 warehouses, and 30 plants. Each algorithm settled into a unique, statistically significant computational time when solving these problems. Sensitivity analyses, supported by non-parametric statistical methods, validate the effectiveness of the proposed heuristic framework.

## 1. Introduction

Supply chain network design is a fundamental strategic problem in operations research. Central to this design is the Multi-Level Facility Location (MFL) problem. In modern logistics, particularly within the retail industry, goods rarely flow directly from producer to consumer. Instead, they traverse a hierarchical network comprising suppliers, plants, warehouses, distribution centers, and retail outlets. The strategic decision of where to locate these facilities and how to assign demand nodes to them dictates the long-term operational efficiency and cost structure of the supply chain.

While MFL models have applications in healthcare and telecommunications, this study focuses specifically on the retail supply chain context. The retail sector is characterized by distinct hierarchical tiers (e.g., suppliers, plants, warehouses, and retailers) and strict single-assignment constraints, often required for administrative simplicity and consistent routing. This context provides a practical testing ground for large-scale optimization, given the massive volume of retail nodes (customers) and the complex upstream dependencies.

Despite its importance, the literature on single-assignment MFL has largely focused on small to medium instances or two-level problems. Solutions to large-scale, realistic problems involving four or five echelons with thousands of customers remain scarce. To address this gap, this paper contributes a comprehensive heuristic framework based on Variable Neighborhood Descent (VND).

The remainder of this article is structured as follows. Section 2 provides a literature review of MFL problems. Section 3 formally defines the problem and presents the Mixed-Integer Linear Programming (MILP) formulation. Section 4 details the VND search variants, including Basic Variable Neighborhood Descent (BVND), Pipe Variable Neighborhood Descent (PVND), Cycle Variable Neighborhood Descent (CVND), and Union Variable Neighborhood Descent (UVND). Sections 5 and 6 present the experimental design, computational results, and statistical validation. Finally, managerial implications and conclusions are discussed.

## 2. Problem Definition

The single-assignment MFL, known as *k*-LFL, is defined as follows. This paper considers 4- and 5-level facilities, i.e., 4-LFL and 5-LFL. In the 4-LFL, we have Levels 1 through 4; in the 5-LFL, Levels 1 through 5. Adapted from the formulation in Ortiz-Astrorquiza et al. (2018), we present below a mathematical formulation of the 5-LFL. In that, we assume that each retailer may be interested only in products from specific plants, and we impose an upper bound on the number of facilities to be opened in Levels 2-5.

$R$: Number of retail stores (Level 1)
$D$: Number of distribution centers (Level 2)
$W$: Number of warehouses (Level 3)
$P$: Number of plants (Level 4)
$S$: Number of suppliers (Level 5)
$r$: A retail store, $r=1,\dots,R$
$d$: A distribution center, $d=1,\dots,D$
$w$: A warehouse, $w=1,\dots,W$
$p$: A plant, $p=1,\dots,P$
$s$: A supplier, $s=1,\dots,S$
$D_r$: The set of distribution centers that can serve a retailer $r$
$P_r$: The set of plants that a retailer $r$ is interested in their products
$W_d$: The set of warehouses that can send products to a distribution center $d$
$P_w$: The set of plants that can send products to a warehouse $w$
$S_p$: The set of suppliers that can send products to a plant $p$
$(s,p,w,d,r)$: A feasible solution for a retailer $r=1,\dots,R$, where: $s \in S_p, p \in P_w, p \in P_r, w \in W_d, d \in D_r$
$M$: A sufficiently large number
$x_l = 1$, if a facility $l$ is opened, 0 otherwise
$v_{d,r} = 1$, if a retailer $r \in R$ is served by a distribution center $d \in D$, 0 otherwise
$z_{(k+1),k,r} = 1$, if a retailer $r \in R$ is served via a route (an edge) from Level $(k+1)$ to Level $k$, 0 otherwise, ($k$=2,…,4)
$u_k$, upper bound for the number of facilities can be selected in Level $k$, ($k$=2,…,5)
$c_{(k+1),k}$, cost of transportation of a bundle of products along a route ($k$+1) to $k$, ($k$=1,…,4)
$F_l$, fixed cost of opening a facility $l$

Min $\sum_{r\in R}\sum_{d\in D_r} c_{dr}v_{dr} + \sum_{r\in R}\sum_{d\in D_r}\sum_{w\in W_d} c_{wd}z_{wdr} + \sum_{r\in R}\sum_{w\in W}\sum_{p\in P_w, p\in P_r} c_{pw}z_{pwr} +$
$\sum_{r\in R}\sum_{p\in P_r}\sum_{s\in S_p} c_{sp}z_{spr} + \sum_{d\in D} F_d x_d + \sum_{w\in W} F_w x_w + \sum_{p\in P} F_p x_p + \sum_{s\in S} F_s x_s$

s.t.

$$\sum_{d\in D_r} v_{dr} = 1, \quad \forall r \in R \tag{1}$$
$$\sum_{w\in W_d} z_{wdr} = v_{dr}, \quad \forall r \in R, d \in D_r \tag{2}$$
$$\sum_{p\in P_w, p\in P_r} z_{pwr} = \sum_{d:w\in W_d} z_{wdr}, \quad \forall r \in R, \mathrm{w} \in W \tag{3}$$
$$\sum_{s\in S_p} z_{spr} = \sum_{w:p\in P_w} z_{pwr}, \quad \forall r \in R, p \in P_r \tag{4}$$
$$\sum_{r\in R:d\in D_r} v_{dr} \le Mx_d, \quad \forall d \in D \tag{5}$$
$$\sum_{r\in R}\sum_{d\in D_r:w\in W_d} z_{wdr} \le Mx_w, \quad \forall w \in W \tag{6}$$
$$\sum_{r\in R:p\in P_r}\sum_{w\in W:p\in P_w} z_{pwr} \le Mx_p, \ \forall p \in P \tag{7}$$
$$\sum_{r\in R}\sum_{p\in P_r:s\in S_p} z_{spr} \le Mx_s, \quad \forall s \in S \tag{8}$$
$$\sum_{l\in Level\ k} x_l \le u_k, \quad \forall k = 2, \dots, 5 \tag{9}$$

All variables binary (10)

The first term in the objective function computes the transportation costs from distribution centers to retailers. The second to fourth terms compute the costs of transporting products along routes in Levels 2-5, for all retailers. Note that retailers are also interested in specific plants. The last four terms calculate the one-time costs of opening facilities in Levels 2-5. Constraint (1) ensures all facilities at Level 1 are served. Constraint (2) ensures that if a distribution center serves a retailer, there is a bundle of products for this retailer transported from a warehouse to that distribution center. Constraints (3-4) are the flow conservation for each retailer along routes between Levels 2 and 5. Here again, notice that constraints (3-4) enforce that a retailer can receive products only from specific plants. Constraints (5-8) require that if a bundle of products is transported to a facility in Levels 2-5, the facility must be opened. Constraint (9) enforces opening and upper bounds for the facilities in Levels 2-5. And Constraint (10) enforces that all variables are binary.

Considering 5-LFL, each retail store $r \in Level\ 1$ must be served via a single-assignment product (a bundle of products), starting from *Level 5* and finally reaching *r*. In that, a bundle of products flows from *Level (k+1)* to *Level k* (for *k*=1,…,4); however, an element in *Level k* can receive products only from a set of elements in facilities in *Level (k+1)* (for *k*=1,…,4).

Furthermore, each retail store $r \in Level\ 1$, besides being interested in receiving a bundle of products from an eligible *Level 2* facility, it may also be interested in receiving products from an eligible set of facilities in *Level 4* (i.e., an eligible facility $p \in Level\ 4$). This may occur in practice, as retailers are often interested in products from specific facilities (e.g., Level 4 plants). A general topology of such a complex network of facilities is illustrated in Figure 1. A schedule for transporting a bundle of products to a retail store $r \in Level\ 1$ is illustrated by path (*s,p,w,d,r*) in Figure 1.

Due to many real factors such as resource limitations and market considerations, we impose upper bounds on the number of facilities selected in each level *k* (*k*=2,...,5). Transporting products through such a complex network incurs some costs. Opening a facility at *Level k* (for *k*=2,…,5*)* involves a one-time fixed cost, while moving a bundle of products from *Level (k+1)* to *Level k* (for *k*=1,…,5) incurs a cost each time it is moved. The objective is to serve all retailers (i.e., elements of the *Level 1* facility) while minimizing total costs. We defined the problem for *k*=5; however, it can similarly be defined for *k*=4.

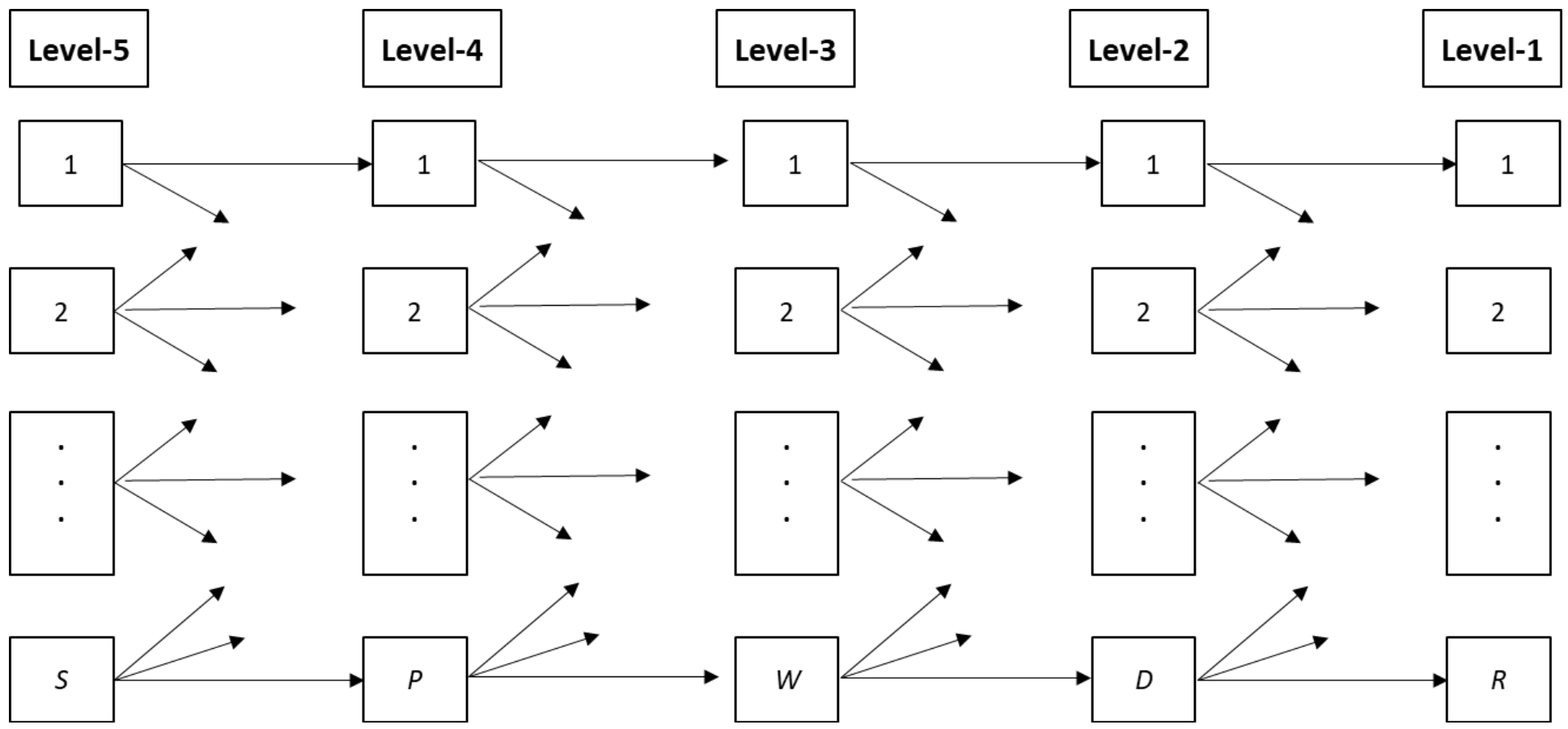

Figure 1. Example of 5-LFL supply chain networks

Due to many real factors such as resource limitations and market considerations, we impose upper bounds on the number of facilities selected in each level *k* (*k*=2,...,5). Transporting products through such a complex network incurs some costs. Opening a facility at *Level k* (for *k*=2,…,5*)* involves a one-time fixed cost, while moving a bundle of products from *Level (k+1)* to *Level k* (for *k*=1,…,5) incurs a cost each time it is moved. The objective is to serve all retailers (i.e., elements of the *Level 1* facility) while minimizing total costs. We defined the problem for *k*=5; however, it can similarly be defined for *k*=4.

This problem has numerous applications in supply chain management. It is similar to the problem discussed in Ortiz-Astorquiza, Contreras, and Laporte (2018, 2019), except that we also consider the interests of retailers with products from a specific set of facilities in *Level 4*. Ortiz-Astorquiza, Contreras, and Laporte (2018) provided a comprehensive review of the MFL. Ortiz-Astorquiza, Contreras, and Laporte (2019) obtained an exact solution to the problem for *k* equal to 2 and 3, solving medium-sized problems. Here, we consider a very large-scale multi-start heuristic based on several variants of the VND meta-heuristic for *k* equal to 4 and 5. In the next section, we present the literature review, followed by a discussion of this paper's contribution.

## 3. Literature Review

This paper focuses on multi-level facility location (MFL) problems. In the literature, MFL is referred to by various terminologies, including multi-echelon facility location, multi-stage facility location, hierarchical facility location, multi-layer facility location, and *k*-level facility location (*k*-LFL). Here, we interchangeably refer to these problems as MFL and *k*-LFL, with *k* equal to 4 and 5. The problem has far-reaching applications in various settings, including tactical and strategic supply chain configuration and transportation planning (Muriel and Simchi-Levi 2003; Melo, Nickel, and Saldanha-da-Gama 2009; Ortiz-Astorquiza, Contreras, and Laporte 2018; Kumar et al. 2020; Janjevic, Merchán, and Winkenbach 2021; Kang, Shen, and Xu 2022; Lin et al. 2022; Wandelt, Wang, and X. Sun 2025). For example, Lin et al. (2022) propose a framework using a discrete choice model and a branch-and-cut algorithm with on-the-fly Benders cuts, based on two analytical separation functions, that appear to efficiently optimize store locations and fulfillment strategies by balancing revenue generation with transportation and fixed operating expenses.

A comprehensive review by Farahani et al. (2014) indicates the model's remarkable adaptability. The model has been used to address public welfare issues, such as determining the optimal locations for new healthcare clinics or positioning emergency medical services to achieve the fastest response times. At the same time, the model can handle complex logistical and infrastructure planning, from organizing municipal solid waste management systems to optimizing the layout of production and distribution networks. Its utility even extends into designing more efficient computer and telecommunication systems. Recent applications in services have been reported in Mogale, Cheikhrouhou, and Tiwari (2020), Majumdar et al. (2023), Kar and Jenamani (2024), Kumar and Kumar (2024), and Ariningsih et al. (2025).

In a recent study of facility location in healthcare, Tsao, Balo, and Lee (2024) emphasize the importance of strategic MFL decisions to enhance supply chain resilience against natural disasters. Janjevic, Merchán, and Winkenbach (2021) focus on multi-facility locations of parcel lockers for last-mile delivery to solve the final, and often most difficult, logistical puzzle of getting a package to a customer.

Although hub-location analysis differs from the hierarchical model mentioned in this paper, it is also considered a multi-echelon approach. It can benefit from our work here, and vice versa. For example, Wandelt, Wang, and Sun (2025) emphasize the multi-echelon location analysis of the e-commerce last-mile delivery system.

Below, we review the literature on single-assignment MFL problems.

The single-assignment uncapacitated multi-level facility location (UMFL) problem generalizes the fundamental single-assignment uncapacitated single-level facility location (UFL) problem. The UFL problem has been extensively researched over many decades (e.g., Laporte, Nickel, and Saldanha da-Gama 2019).

An early study of the MFL problem was conducted by Kaufman, Eede, and Hansen (1977), who introduced the two-level so-called warehouse and plant location problem. Ortiz-Astorquiza, Contreras, and Laporte (2018) and Kumar et al. (2020) noted that most research on MFL has focused on two- or three-level cases, which is also discussed in Malik, Contreras, and Vidyarthi (2022), Sluijk et al. (2023), Fokouop et al. (2024), Gendron, Khuong, and Semet (2024), and Marianov and Eiselt (2024). We refer to a comprehensive review of general MFL problems by Ortiz-Astorquiza, Contreras, and Laporte (2018), including single-assignment cases.

Chardaire, Lutton, and Sutter (1999) formulated a telecommunication problem as a single-assignment two-level model and provided upper and lower bounds for the problem. Previously, Tragantalerngsak, Holt, and Rönnqvist (1997) also formulated a transportation problem as a single-assignment two-echelon model and applied a Lagrangian heuristic. Gendron, Khuong, and Semet (2015) implemented a multilayer variable neighborhood search (VNS) to solve a two-level, single-assignment facility location problem. Gendron, Khuong, and Semet (2016) presented a Lagrangian branch-and-bound approach for the optimal solution to the two-level uncapacitated single-assignment problem.

Hammami, Frein, and Bahli (2017) examined the impact of lead time on a single-assignment, multi-echelon facility location in the supply chain. They concluded that manufacturing and distribution sites should be located near the demand zone (i.e., retailers) and local suppliers should be selected, despite their higher costs.

Ortiz-Astorquiza, Contreras, and Laporte (2018, 2019) studied the UMFL with single-assignment constraints. Path-based and arc-based integer programming approaches are given. An optimal solution approach based on Benders' reformulation is presented, capable of solving medium-sized two- and three-level problems.

Ramshani et al. (2019) considered the single-assignment, two-level, uncapacitated facility location problem with uncertain disruptions. The authors developed two mathematical models and a tabu search as a solution approach.

Myung and Yu (2020) analyzed a freight transportation model with a multi-echelon structure that bundles products. They mentioned that transporting product bundles can be viewed as a special case of single-assignment; however, additional effort and costs are associated with bundling and unbundling products. The authors developed a heuristic based on the network flow algorithm.

De Oliveira et al. (2021) examined how well a specific "divide-and-conquer" algorithm works for a dynamic facility location problem. Their study examined a two-level system without capacity limits, testing the algorithm's performance when customers must be assigned to a single facility and when they can be assigned to multiple facilities.

Wang et al. (2023) addressed a single-item, multi-echelon location inventory and provided a $\varepsilon -$ optimal approach using Lagrangian relaxation. Although the hub facility location is hierarchical and shares a significant similarity with the MFL considered in the current study, they are not identical. Our approaches can benefit hub location analysis and vice versa. For a comprehensive review of articles from an air transportation perspective, refer to O'Kelly, Sun, and Wandelt (2025), which includes a review of single-assignment models.

**Mathematical Formulation**

We formulate the single-assignment $k$-Level Facility Location ($k$-LFL) problem. The supply chain consists of $K$ levels of facilities. Let Level 1 represent the customers (retailers) and Level $K$ represent the highest echelon (e.g., suppliers). The goal is to minimize total fixed opening costs and transportation costs while ensuring every node at Level $l$ is assigned to exactly one open node at Level $l + 1$.

**Sets and Indices:**

- $L = 1, \dots, K$: Set of facility levels.
- $N_l$: Set of potential facility locations at the level $l \in L$.
- $i, j$: Indices for facility locations.

**Parameters:**

- $F_j$: Fixed cost of opening a facility $j$ at level $l$ (where $l > 1$).
- $C_{ij}$: Cost of assigning a node $i \in N_l$ to a facility $j \in N_{l+1}$. This represents the cost of moving a product bundle between levels.
- $U_l$: Maximum number of facilities allowed to open at the level $l$.

**Decision Variables:**

- $y_j \in 0,1$: Equal to 1 if facility $j$ is opened; 0 otherwise. (Note: Level 1 nodes are always "open" as demand points).
- $x_{ij} \in 0,1$: Equal to 1 if node $i \in N_l$ is assigned to the node $j \in N_{l+1}$; 0 otherwise.

**Model:**

$$\text{Minimize } Z = \sum_{l=2}^{K} \sum_{j \in N_l} F_j\, y_j + \sum_{l=1}^{K-1} \sum_{i \in N_l} \sum_{j \in N_{l+1}} C_{ij} x_{ij}$$

Subject to:

$$\sum_{j \in N_{l+1}} x_{ij} = 1 \forall l \in 1, \dots, K-1, \forall i \in N_l$$

(Single assignment constraint: each lower-level node assigns to exactly one upper-level node).

$$x_{ij} \le y_j \forall l \in 1, \dots, K-1, \forall i \in N_l, \forall j \in N_{l+1}$$

(Connectivity constraint: nodes can only be assigned to open facilities).

$$\sum_{j \in N_l} y_j \le U_l \forall l \in 2, \dots, K$$

(Cardinality constraint: upper bound on open facilities per level).

$$x_{ij}, y_j \in 0{,}1$$

This formulation captures the 4-LFL ($K = 4$) and 5-LFL ($K = 5$) problems discussed in this study. The problem is NP-hard, motivating the use of the VND heuristics described in Section 5.

## 4. Contribution of This Paper

As the literature suggests, for single-assignment MFL, most algorithms are based on two-level (rarely three-level) facility location problems. Even in these cases, small to medium-sized problems are being considered. Real-world problems are large-scale, often involving more layers (e.g., supply chain design problems), and within each layer, potentially more facilities are considered (Zandi Atashbar, Labadie, and Prins 2018). Note that even in the simplest case, the ULFL is NP-hard (nondeterministic polynomial time hard). MFL problems are more complex and require sophisticated algorithms for large-scale cases. The problems considered in this study are not only larger in scale, but they also involve additional constraints related to the literature, as explained above. To address these shortcomings, we consider the VND meta-heuristic, a variant of VNS, for large-scale 4- and 5-LFL problems. Although there are many variants of VND, the most representative are BVND, PVND, CVND, and UVND (Duarte et al. 2018). The detailed explanation of these four variants for 5-LFL is explained in the next section. Again, tuning the algorithms for 4-LFL is straightforward.

A brief overview of recent applications of the four variants of VND (BVND, PVND, CVND, and UVND) algorithms is provided below. As mentioned earlier, Duarte et al. (2018) provided an excellent discussion of these VND variants and their applications. Liu et al. (2021) also proposed a hybrid VND genetic-particle swarm algorithm for the flexible job shop. The VND is integrated into the standard genetic algorithm (GA) framework as an improvement process. Mjirda et al. (2017) employed a sequential VND for the traveling salesman problem. The authors used the four variants we also consider here as the moves in their procedures. Daquin et al. (2021) applied two VND variants, BVND and UVND, as improvements over general VNS and used them to cross-dock truck assignments. They employed the Best Improvement (BI) and First Improvement (FI) processes in the algorithms. Siew, Sze, and Goh (2025) considered VND variants, including the four variants we consider here, as an improvement process in the Whale Optimization Algorithm (WOA). Surprisingly, variants of VNS are absent from applications to MFL problems.

## 5. Four Variants of VND Meta-heuristic for MFL

An optimization problem may be defined by a feasible solution *X*, and an objective function $f: X \rightarrow Real$ where *Real* is the set of real numbers. The problem is to find a solution, $x* \in X$ that optimizes (maximizes or minimizes, depending on the problem) the function *f*. Without loss of generality, we focus on minimization, noting that maximization is equivalent to minimizing *-f*.

Members of the set *X* may be defined differently depending on the problem considered. Often, $x \in X$ is a binary vector, or a vector of real numbers.

Below, we provide several definitions for 5-LFL; similar definitions may also apply to 4-LFL. In many cases, the set of solutions is defined on a network, as is the case here. Earlier, we explained that a feasible solution (schedule) for a customer $r \in Level\ 1$ is defined by a path (*s,p,w,d,r*) in the network shown in Figure 1. We show this individual solution by $x_r = (s, p, w, d, r)$. Generally, many feasible solutions for a retailer *r* create a set of solutions denoted by $X_r$. Now, we can define the set of feasible solutions for the problem as $X = X_r$ for *r*=1, …, *R*.

Given a feasible solution, $x_r \in X_r$ for a retail store $r \in Level\ 1$, let $x'_r \in X_r$ when *k* (for *k*=1,…,4) eligible elements from the set {*s,p,w,d*} are changed. The set of all such $x'_r \in X_r$ is called the set of *k*-neighborhoods of $x_r$, denoted by $N_k(x_r)$. For a retail store *r*, $N_k(x_r)$ includes 4 (one-element change), 6 (two-element change), 4 (three-element change), and 1 (four-element change), also called neighborhood types, for *k*=1, 2, 3, and 4, respectively, for a total of 15 neighborhood types. Furthermore, $\mathbb{N}(k) = \cup_{r=1}^{R} N_k(x_r)$ is referred to as the set of *k*-neighborhood structures for the problem. Let $x_{loc} \in \mathbb{N}(k)$ be such that $f(x_{loc}) \le f(x)$ for all $x \in \mathbb{N}(k)$, $x_{loc}$ gives us the best possible outcome. Let $\mathbb{N} = \cup_{k=1}^{max\ \ k} N(k)$, and $x* \in \mathbb{N}$ be such that $f(x^*) \le f(x)$ for all $x \in \mathbb{N}$, x* is a global optimal solution to the problem.

Over the last several decades, researchers have formulated a wide range of real-world problems as optimization problems and solved them using various methods. Since these problems are NP-hard and solving realistic, large-scale problems is difficult, researchers often use heuristics, specifically metaheuristics. Duarte et al. (2018) noted that over the last several decades, more than 50 variants of meta-heuristics have been developed. In general, such heuristics can be categorized into Adaptive Heuristics (AH), Large Neighborhood Search (LNS), Adaptive Large Neighborhood Search (ALNS), and various combinations of these methods.

AH is a heuristic that modifies its configuration as the search progresses, thereby changing its behavior (Gouda and Herman 1991; Sevaux, Sörensen, and Pillay 2018; Yaakoubi and Dimitrakopoulos 2025). First introduced by Shaw (1998), LNS is a metaheuristic built on a surprisingly straightforward principle: to find a better solution, it may be necessary to selectively break an existing good one. The core idea isn't to construct a solution from scratch, but to operate through a cyclical process of partial destruction and subsequent reconstruction. At each search step, the algorithm typically pairs one of each. It first takes a destructive approach to break down the current solution, resulting in a partial, incomplete version. Then, using a constructive method, it rebuilds from that point, ideally charting a path to an improved overall outcome. Often, the selection of constructive and destructive processes is implemented probabilistically. ALNS, initially developed by Ropke and Pisinger (2006), builds on LNS by probabilistically selecting a pair of constructive and destructive processes based on information from previous runs.

In general, the basic idea of AH involves dynamic neighborhood changes to reach better solutions when a simpler neighborhood is stuck in a local optimum. This dynamism expands the search to diverse regions of the search space to find better solutions. Refer to Sevaux, Sörensen, and Pillay (2018) for a general discussion of AH models, Ahuja et al. (2002) for a survey of very large-scale neighborhood search algorithms.

A popular adaptive meta-heuristic is the VNS, proposed by Mladenović and Hansen (1997). The VNS relies on several local search processes across several neighborhood structures. The General VNS is a type of ALNS and, thus, often depends on the pair of *constructive* and *destructive* processes in a solution. However, several variants of the VNS have been proposed in the literature, including VND, Reduced VNS (RVNS), Basic VNS (BVNS), General VNS (GVNS), Skewed VNS (SVNS), and Variable Neighborhood Decomposition Search (VNDS); refer to Duarte et al. (2018) for a basic discussion of these methods. Among these variants, VND is one of the most widely used by researchers; it is effective and straightforward to implement and does not rely on constructive or destructive phases. Over the years, many variants of VND have been proposed; see, for example, Mjirda et al. (2017), Duarte et al. (2018), Daquin et al. (2021), and Siew, Sze, and Goh (2025).

In the VND algorithms, a list of neighborhood structures is provided sequentially, typically to increase sophistication (from the smallest to the largest). Within these variants, Basic VND (BVND), Pipe VND (PVND), Cyclic VND (CVND), and Union VND (UVND) are the most representative, as noted by Duarte et al. (2018). These methods differ in the order and depth at which neighborhood structures are applied. This means that when an improvement in a neighborhood is detected, it is important to determine which neighborhood to explore next and how deeply to investigate it.

It is believed that the larger the neighborhood, the better the quality of the locally optimal solutions. However, using a larger neighborhood incurs high CPU time costs (Ahuja et al. 2002). Thus, when applied to large-scale problems, it is critical to balance CPU time and to explore a larger neighborhood to find better solutions. A general description of these algorithms is given below, focusing on where to explore next after an improvement is detected in the process (Duarte et al. 2018).

After finding a way to improve the solution:

- BVND – This strategy returns to the beginning and reapplies its first search method on the newly improved solution.
- PVND – This strategy uses the same search method that just succeeded.
- CVND – This strategy proceeds to the next search method in its predefined list.
- UVND – This strategy treats all search methods as one big "toolbox." After finding an improvement with one tool, it continues searching using any tools from its entire collection, without a strict order.

As mentioned, the neighborhood structures are presented in order of sophistication. Thus, the BVND process is fully explored each time the least sophisticated (simplest) neighborhood structure is reached. In the PVND, each neighborhood structure is fully explored before proceeding to the next structure. In the CVND, however, as soon as an improvement is detected, the search continues to the next structure. In the case of UVND, all structures are combined into a single large structure and appropriately explored.

An important factor to consider when implementing a local search within a specific neighborhood is the order in which the search process is executed. This is especially important when dealing with very large-scale problems. In these cases, the order of implementation is effective in reaching actionable results. This provides an opportunity to explore a more diverse area of the solution space. It has been experimentally demonstrated that a random order of implementations yields significantly better results than always choosing a specific order. Choosing a random order each time to explore the neighborhood structure can be time-consuming. However, selecting a sequence each time can significantly reduce this time, as discussed in Wang and Alidaee (2023). We explore this factor in more detail later in the paper, when we describe the algorithms' pseudocode.

In the following subsection, we provide the pseudocode for the 5-LFL algorithm, which can be easily adapted for 4-LFL.

### *5.1 Pseudocode of the Four Algorithms for 5-LFL*

Considering 5-LFL, given a schedule $x_r = (s, p, w, d, r)$ for a retail store *r*, with all neighborhood structures $N_k(x_r)$ for *k*=1,…,4, and *k* possible changes among *s, p, w*, and *d*. We also refer to this as a *k*-flip move (or exchange). Thus, for *k*=1, 2, 3, and 4, the number of elements in $N_k(x_r)$, respectively, is *n*=4, 6, 4, and 1, totaling 15 neighborhood moves.

As previously mentioned, a critical issue in designing heuristics is the choice of neighborhood structure and the order in which they are implemented. The order of implementation for a local search within each neighborhood structure is also an important factor. Furthermore, in optimization problems, neighborhood structures are often defined based on the problem's mathematical programming formulation. However, a graphical neighborhood structure is considered in this context.

In the case of 4-LFL, we have three neighborhood structures $\mathbb{N}(k)$ for *k*=1, 2, and 3. For $\mathbb{N}(k)$ and *k*=1 and 2 each include three neighborhood types, while for k=3 there is only one, for a total of seven neighborhood types, as illustrated in Figure 2. In the case of 5-LFL, we have four neighborhood structures $\mathbb{N}(k)$ for *k*=1, 2, 3, and 4. For $\mathbb{N}(k)$ and *k*=1 and k =3, each include four neighborhood types, while for k=2, it includes six neighborhood types, and for *k*=4, it includes one, for a total of 15 neighborhood types, as illustrated in Figure 3. These figures are shown in Appendix A.

To implement the *k*-flip move processes, for a given value of *k*, several important issues should be considered:

(a) Which retail stores, *r,* should be considered each time for a possible *k*-flip move implementation? We use a sequence *Lr(1), …, Lr(R)* of *R* numbers to select the next retail store.
(b) Given a schedule for a retail store, *r*, which combination of *k* elements from the set $N_k(x_r)$ should be considered for *k*-flip moves? For this, we use a sequence, *q(1),…,q(n),* for example, for *k*=2, we use a sequence of *n*=6 numbers *q(1),…,q(6)*, indicating an order of six elements in $N_2(x_r)$. Thus, for *k*=2, each element $q(.) \in N_k(x_r)$ is a pair of nodes.
(c) Given a schedule for a retail store, *r*, a neighborhood structure, $N_k(x_r)$ for some value of *k*, and an element $q(.) \in N_k(x_r)$, which element of *q*(.) should be considered next for a *k*-flip move? For example, for (*s,p,w,d,r*) and *k*=2, let $q(.) = d, p$. For this, we consider sequences *Ld(1),…,Ld(D)*, and *Lp(1),…,Lp(P)*, respective, *D* and *P* numbers. Thus, along these two sequences, we flip two nodes $d' \neq d$ and $p' \neq p$ for a possible 2-flip move. Similarly, we can implement *k*-flip moves for different values of *k*.

Algorithm-0 is a basic *k*-flip local search strategy appropriately used in other algorithms.

**Algorithm-0:** Simple *k-flip* Neighborhood Search Process (5-LFL)
**Initialization:** Set of numbers *R, D, W, P,* and *S*.
A value for $k$ (1,…,4), and the set of neighborhood structures $\mathbb{N}(k)$. A feasible schedule (*s,p,w,d,r*) for each retail store *r*=1,…,*R*.
*Improvement*=True
**While** (*Improvement*) **Do**
  *Improvement*=False
  Randomly select an order of numbers 1,…,*n, i.e., q(1),…,q(n),* of elements in $\mathbb{N}(k)$
  **For (*h=q(1),…,q(n)*)**
    Randomly select sequences, *Lr, Ls, Lp, Lw,* and *Ld,* of numbers *R, S, P, W*, and *D*, respectively.
    If a *h-flip* is improving along appropriate sequences *Lr, Ls, Lp, Lw,* and *Ld,* implement the move, and set *improvement*=True
  **End For**
  Update the best-known solution
**End While**

Algorithm-0 exhaustively implements the process of *k*-flips until no more moves are possible. In Algorithm-0-k, however, the algorithm returns when an improvement is detected.

***Algorithm-0-k(.)*** for *k*=2,3,4, is the same as Algorithm-0 except that as soon as an improvement is found for any of the *h* values, RETURN the result.

**Algorithm-1 (Multi-start *k*-flip):** Multi-start Neighborhood Search Process (5-FLP)
**Initialization:** Set of neighborhood structures, $\mathbb{N}(k)$, for a value of *k*=1,…,4. Set of integer numbers *R, D, W, P,* and *S.* Max-Local *(number of multiple starts), j*=1
**While** (*j<=Max-Local*) **Do**
  Randomly select sequences, *Lr, Ls, Lp, Lw,* and *Ld*
  Call **Algorithm-0** for *k*
  *j=j+1*
  Keep track of the best solution found throughout
**End While**

**Algorithm-2:** BVND Multi-start Neighborhood Search Process (5-LFL)
**Initialization:** Set of neighborhood structures, $\mathbb{N}(k)$, (for *k*=1,…,4)
Call *Algorithm-1(.)* with *k*=1
*Improvement*=True
**While** (*Improvement*) **Do**
  *Improvement*=False
  **Step 1.** Call *Algorithm-0(.)* with *k*=1.
  **For** (*k*=2,…,4) **Do**
    Call *Algorithm-0-k(.)* with *k*, if an improvement is detected, implement the change, *improvement*=True, and **go to Step 1**, otherwise **continue**
  **End For**
  Update the best-known solution
**End While**

Note that in **Step 1** of Algorithm-2, the 1-flip Algorithm-0(.) is exhaustively implemented; however, the same algorithm for *k*=2, 3, and 4 returns to **Step 1** as soon as an improvement is detected, i.e., Algorithm-0-k(.). Also, note that Algorithm-1(.) with *k*=1 is the multi-start 1-flip local search strategy. The result of this algorithm is used as a starting solution for the BVND, PVND, CVND, and UVND algorithms.

The choice of different sequences in the algorithms is crucial. Using different sequences allows diversification into a large area of the solution space as the search progresses. It is also easy to start different solutions in the multi-start strategy. Any method for selecting new sequences each time is acceptable. However, a clever implementation of the sequence selection process can significantly reduce CPU time. Two such innovative approaches are adapted from applications in sequencing problems, such as the Traveling Salesman Problem (TSP). One is based on the so-called *l*-Opt local search strategy applied to the TSP and many other sequencing problems (see Alidaee and Wang 2017). The other is the Random-Key strategy (see Wang and Alidaee 2023), adapted from the Random-Key application in the TSP (Bean 1994). There are both advantages and disadvantages to using these two strategies. Using the *l*-Opt strategy may be more time-consuming, but it uses less memory. However, the opposite is true for the Random-Keys application. To address very large-scale problems, we employed the *l*-Opt strategy to minimize space usage. Refer to Alidaee and Wang (2017) for a detailed implementation of the *l*-Opt strategy, and Wang and Alidaee (2023) for details on the Random-Key strategy supplication. Note that the value of $r$ in the *l*-Opt strategy is effective in results when solving problems. We used a limited 4-Opt strategy in our implementation, adapted from Glover (1996) for TSP applications.

As explained earlier, there are different ways to define neighborhood structures. Often, they are determined based on mathematical programming formulations of the problem. However, this study uses graphical structures, as shown in Figure 2 for 4-LFL and Figure 3 for 5-LFL.

We previously noted that variations of the VND methods differ in the order and depth in which the neighborhood structures are implemented. This can create many variants, some of which could be very time-consuming. However, four popular variants are BVND, PVND, CVND, and UVND. Even within each of these variants, there are many ways to implement the process. After applying several approaches (see the computational experiments section), it was found that implementing a multi-start strategy, similar to Algorithm 1, significantly reduces CPU time. At the same time, the quality of the results of Algorithms 2 through 5 also remains high.

A *k-flip local search* can be exhaustively implemented for each neighborhood structure $\mathbb{N}(k)$. Generally, implementing local searches for larger values of $k$ is more time-consuming. However, it is also believed that larger values of $k$ can diversify to broader areas of the solution space and possibly create better solutions. Thus, it is important to balance CPU time and achieve better solutions. Here, we used a multi-start strategy, Algorithm-1, which incorporates the 1-flip strategy. The results of this algorithm serve as the starting solution in Algorithms 2 through 5.

Algorithm 2 is written for BVND; however, it can easily be tuned for other VND variants. Algorithms 3, 4, and 5 are designed for PVND, CVND, and UVND. Note that in PVND, before the search moves from one neighborhood structure to another, an exhaustive 1-flip local search is completed. However, in CVND, the search is explored in the next neighborhood structure after each local improvement is detected. Also, note that in the UVND, it is important to specify the order of the local improvement process within the *big* neighborhood structures. The union of all neighborhoods includes possible moves of 1, 2, 3, or 4 moves. Here, we randomly permuted the order to assess the improvement in moves within the *big* neighborhood structure. The UVND has similarities with CVND. In CVND, we implement neighborhood structures in a specific order; however, in UVND, we treat all neighborhood structures equally and randomly apply the improvement process. Also, note that most researchers use the so-called Best Improvement (BI) or First Improvement (FI) process in UVND (Daquin et al., 2021); however, we use the Next Improvement process here.

**Algorithm-3:** PVND Multi-start Neighborhood Search Process (5-LFL)
**Initialization:** Set of neighborhood structures, $\mathbb{N}(k)$, (for $k$=1,…,4)
Call *Algorithm-1(.)* with *k=1*
*Improvement*=True
**While** (*Improvement*) **Do**
  *Improvement*=False
  **For** ($k$=1,…,4) **Do**
    Call *Algorithm-0(.)* with $k$, if an improvement is detected, implement the change.
    *Improvement*=True
  **End For**
  Update the best-known solution
**End While**

Note that Algorithm-3 does not return immediately after finding an improvement in Algorithm-0(.) with $k$. It exhaustively continues the search in the same neighborhood structure, then returns and restarts Algorithm-0(.) with $k+1$.

**Algorithm-4:** CVND Multi-start Neighborhood Search Process (5-LFL)
**Initialization:** Set of neighborhood structures, $\mathbb{N}(k)$, (for $k$=1,…,4)
Call *Algorithm-1(.)* with *k=1*
*Improvement*=True
**While** (*Improvement*) **Do**
  *Improvement*=False
  **For** ($k$=1,…,4) **Do**
    Call *Algorithm-0-k(.)* with $k$, if an improvement is detected, implement the change
    *Improvement*=True,
  **End For**
  Update the best-known solution
**End While**

Note that, in CVND, Algorithm-0-k(.) is used in the inner loop, as each time a local search is detected, we return to the next neighborhood. Additionally, the inner loop is consistently implemented with the same order: $k$ = 1, 2, 3, 4. You may obtain different results for each order if other orders are used. However, CPU time will increase.

**Algorithm 5:** UVND Multi-start Neighborhood Search Process, (5-LFL)
**Initialization:** Set of neighborhood structures, $\mathbb{N}(big) = \mathbb{N}(1) \cup \mathbb{N}(2) \cup \mathbb{N}(3) \cup \mathbb{N}(4)$
Call *Algorithm-1(.)* with *k=1*
*Improvement*=True
**While** (*Improvement*) **Do**
  *Improvement*=False
  Randomly select an order of four neighborhood structures, *M(j)*, $j$=1,…,4, in $\mathbb{N}(big)$
  **For** (*k=M(1),…,M(4)*) **Do**
    Select an order of numbers 1,…,$n$, (i.e., *q(1),…,q(n)*), where $n$ is the number of elements in $\mathbb{N}(k)$.
    **For** (*h=q(1),…,q(n)*) **Do**

Call *Algorithm-0-k(.)* with *k*, if an improvement is detected, implement the change, *Improvement*=True, and **go to Start a new neighborhood**, otherwise **continue**
**End For**
**Start a new neighborhood**
**End For**
Update the best-known solution
**End While**

In UVND, the loop running via $k$ goes through all neighborhood structures, and the loop running via $h$ goes through each element of $\mathbb{N}(k)$. As soon as an improvement is detected, it proceeds to the next neighborhood. It should be clear that the orders used in the implementation processes can significantly impact the outcomes in UVND and other algorithms, particularly in terms of solution values and CPU time. Thus, there are many ways to create an effective variant of VND.

## 6. Experimental Design and Results

### *6.1 Data Generation*

No benchmarks are available for the problems considered in this paper. The only benchmark that shares some characteristics with our problems is provided by Ortiz-Astorquiza, Contreras, and Laporte (2019). However, there are only two with four-level facilities. These two problems also lack some data that cannot be used in our computational experiment. Thus, we randomly generated problem instances and solved them using algorithms. All algorithms were implemented in Fortran and executed in order on a cluster with Intel Haswell Xeon processors.

We generated problems of varying sizes with different matrix densities and fixed cost levels. Table B1 in Appendix B shows the parameters used to generate the data. For each problem size (High, Medium, and Low densities) and each fixed cost (Large, Medium, and Small), we generated three instances and solved them using the B-VND, P-VND, C-VND, and U-VND algorithms. Results for 4-LFL are shown in Table 1, and for 5-LFL are shown in Table 2. Note that problem IDs for 4-LFL are shown by (R-D-W-P-Density-Fixed Cost-#). For example, (2000-150-50-30-Hdens-LgFx-1) means R=2000, D=150, W=50, and P=30, with High Density, Large Fixed Cost, and problem # 1. Similarly, problem IDs for 5-LFL are shown by (R-D-W-P-S-Density-Fixed Cost-#). For example, (2000-150-50-50-100-Hdens-LgFx-1) means R=2000, D=150, W=50, P=50, S=100, High Density, Large Fixed Cost, and problem number 1. Three instances of each problem were generated and solved using B-VND, P-VND, C-VND, and U-VND. The objective function and the CPU time required to reach the best solution are given for each case. All instances are available at https://figshare.com/s/e5a4bcd1d35612230f92.

[Insert Tables 1 and 2 here]

### *6.2 Computational Results*

As discussed earlier, each VND variant can be implemented in many different ways. However, it is important to strike a balance between CPU time and solution quality. This is especially important when dealing with large-scale problems. The multi-start strategy provides an opportunity to initiate a new solution each time, resulting in distinct end solutions and potentially higher-quality solutions. To balance CPU time with the final solution, we applied a multi-start strategy on Algorithm-1 with $k$=1, a speedy process. Then, we used the best result of this process as a starting solution for each variant of the VND process. Tables 1 and 2 show the results. To evaluate the solution quality of the VND heuristics, we solve seven 5-LFL instances with R=100 and R=200 using the state-of-the-art Gurobi 13.0 solver to obtain high-quality solutions within a 7200-second time limit, and compare the heuristics' solutions with a 60-second time limit. The details of instances and solution comparison are provided in Appendix C. For instance, with R> 500, the Gurobi 13.0 solver is unable to obtain a feasible solution within the available CPU memory limit and the 18000-second time limit. Only solutions from VND variants will be reported.

### *6.3 Sensitivity Analyses of Algorithms*

To determine which of our algorithms performs better across different tasks, we employed rank-based statistical methods, which are well-suited for this type of evaluation. Instead of getting overwhelmed by precise numerical outcomes, these tests focus on the relative ordering of algorithm performance across problems. Of course, care was taken to ensure these comparisons were meaningful. This includes running sufficient experiments and being thoughtful about our significance thresholds. It is especially important to account for multiple comparisons, so we report not only the raw significance values (p-values) but also, where appropriate, effect sizes to provide a clearer picture of the results.

Our first step was to get a bird's-eye view of the overall performance. For this, we used Friedman's test, a nonparametric method that serves a similar purpose to repeated-measures ANOVA. It works by ranking algorithms on each dataset and then checking whether the average ranks differ too much to be explained by random chance. A significant result from Friedman's test suggests that at least one algorithm behaves differently but does not specify which ones.

When the Friedman test indicated a meaningful difference, we needed to dig deeper to identify the specific pairs of algorithms that were outperforming others. The Nemenyi test is designed for this situation, comparing all algorithms against each other while carefully controlling the error rate that can arise from so many simultaneous comparisons. In other cases, where we were more interested in direct head-to-head contests between two specific algorithms, the Wilcoxon signed-rank test was a better fit. This approach is useful because it considers the direction and magnitude of the performance differences. When running multiple Wilcoxon tests, we adjusted our significance levels using methods such as the Bonferroni correction to avoid drawing false conclusions from multiple tests.

Key statistics for the primary omnibus test are reported in Tables 3, 4, 5, and 6, including the test statistics, degrees of freedom, and the p-value, thereby clarifying the overall findings. Next, we specify which pairs of algorithms showed statistically significant differences using the follow-up methods we applied, making it easy to identify where meaningful distinctions lie.

**Table 3.** Results of Friedman's Test, Nemenyi's Post-Hoc Test, and Pairwise Wilcoxon Test on the OFV of 4-LFL instances.

| **Friedman's Test for BVND, PVND, CVND, and UVND** |
|---|
| Number of complete experiments (blocks/rows): 162 |
| Number of algorithms (groups/columns): 4 |
| Friedman chi-squared statistics: 397.855 |
| P-value: 6.454E-86 |

Note: The Friedman test was conducted to determine if there were any statistically significant differences in the median Objective Function Values (OFVs) achieved by the four algorithms (BVND, PVND, CVND, and UVND) across the 162 experiments. With a chi-squared statistic of 397.855 and an extremely low p-value (6.454E-86), well below the 0.05 significance level, the null hypothesis that all algorithms perform equally well (i.e., have similar median OFVs) was rejected. This indicates that at least one algorithm's typical OFV performance is statistically different from the others, necessitating post-hoc analysis to identify specific pairwise differences.

**--- Post-Hoc Analysis ---**

**Nemenyi's Post-Hoc Test Results (pairwise p-values):**

| | BVND | PVND | CVND | UVND |
|---|---|---|---|---|
| BVND | | 0 | 0 | 0.931 |
| PVND | 0 | | 1.02E-11 | 0 |
| CVND | 0 | 1.02E-11 | | 0 |
| UVND | 0.931 | 0 | 0 | |

Note: Following the significant Friedman test, Nemenyi's post-hoc test was used for pairwise comparisons of algorithm performance based on their OFVs, with a significance level of 0.05. The results indicate that the OFVs achieved by the algorithms BVND and UVND are not statistically distinguishable (p = 0.931). However, all other pairwise comparisons yielded p-values < 0.05, indicating statistically significant differences in their OFV performance. Specifically, PVND and CVND each have OFV distributions that differ significantly from each other, and both differ significantly from the BVND/UVND pair, suggesting three distinct tiers of OFV performance among the algorithms.

**Pairwise Wilcoxon Test with Bonferroni Correction (*p*-values)**

| | BVND | PVND | CVND | UVND |
|---|---|---|---|---|
| BVND | | 1.48E-27 | 1.29E-26 | 1.00 |
| PVND | 1.48E-27 | | 3.39E-26 | 2.95E-26 |
| CVND | 1.29E-26 | 3.39E-26 | | 2.05E-25 |
| UVND | 1.00 | 2.95E-26 | 2.05E-25 | |

Note: The pairwise Wilcoxon signed-rank tests, with a Bonferroni correction to maintain a family-wise error rate of 0.05, were also conducted to compare OFV performance between pairs of algorithms. These results corroborated the findings from Nemenyi's test: no statistically significant difference in the OFVs produced by BVND and UVND (p = 1.00). Conversely, all other pairwise comparisons (BVND vs. PVND, BVND vs. CVND, PVND vs. CVND, PVND vs. UVND, and CVND vs. UVND) showed highly statistically significant differences in their OFV distributions (all adjusted p-values < 0.05, many << 0.001). This reinforces the conclusion that PVND and CVND differ, and that the BVND/UVND pair is statistically similar in terms of the OFVs they yield.

**Table 4.** Results of Friedman's Test, Nemenyi's Post-Hoc Test, and Pairwise Wilcoxon Test on the computing time of 4-LFL instances.

**Friedman's Test for BVND, PVND, CVND, and UVND**

Number of complete experiments (blocks/rows): 162
Number of algorithms (groups/columns): 4
Friedman chi-squared statistic: 332.556
P-value: 8.926E-72

Note: The Friedman test was employed to assess whether there were statistically significant differences in the median computing times among the four algorithms (BVND, PVND, CVND, and UVND) across the 162 experiments. The test yielded a chi-squared statistic of 332.556 and a p-value of 8.926E-72. Since this p-value is less than the pre-defined significance level ($\alpha = 0.05$), the **null hypothesis ($H_0$)**, stating that all algorithms have the same median computing time, is **rejected**. This indicates strong evidence that at least one algorithm's computing time distribution differs significantly from the others, warranting post-hoc analysis to identify specific pairwise differences.

**--- Post-Hoc Analysis ---**

**Nemenyi's Post-Hoc Test Results (pairwise p-values):**

| | BVND | PVND | CVND | UVND |
|---|---|---|---|---|
| BVND | | 0 | 0 | 1.14E-10 |
| PVND | 0 | | 0.678 | 0 |
| CVND | 0 | 0.678 | | 5.22E-15 |
| UVND | 1.14E-10 | 0 | 5.22E-15 | |

Note: Nemenyi's post-hoc test was conducted to perform pairwise comparisons of algorithm computing times following the significant Friedman test, using α = 0.05. The results show that the computing times for PVND and CVND algorithms are not statistically distinguishable (p = 0.678). However, all other pairwise comparisons yielded p-values well below 0.05 (e.g., BVND vs. PVND, p=0.00; BVND vs. UVND, p=1.14E-10; PVND vs. UVND, p=0.00; and CVND vs. UVND, p=5.22E-15). This suggests that algorithms BVND and UVND each have significantly different computing times than those of all other algorithms. At the same time, PVND and CVND form a group with similar computing times.

**Pairwise Wilcoxon Test with Bonferroni Correction (p-values)**

| | BVND | PVND | CVND | UVND |
|---|---|---|---|---|
| BVND | | 4.74E-27 | 3.34E-27 | 2.68E-20 |
| PVND | 4.74E-27 | | 1.13E-02 | 5.75E-24 |
| CVND | 3.34E-27 | 1.13E-02 | | 3.31E-10 |
| UVND | 2.68E-20 | 5.75E-24 | 3.31E-10 | |

Note: A more detailed analysis, carefully adjusted for multiple comparisons to ensure reliability, examined the differences in computing time between every pair of algorithms (BVND vs. PVND, BVND vs. CVND, etc.). This analysis found a real, meaningful difference in computing time between *all* pairs of algorithms. This means BVND's computing time differed from PVND's, CVND's, and UVND's. Similarly, PVND's computing time was different from CVND's (though this difference, while real, was less strongly evident than others), and UVND's and CVND's computing time was also different from UVND's.

The four algorithms exhibit distinct performance profiles in terms of solution quality (OFV) and computational speed on 4-LFL instances. Algorithms BVND and UVND achieve statistically similar solution qualities, but their computing times differ significantly from each other and from those of the other two algorithms. In contrast, PVND and CVND each deliver unique solution qualities, distinct from each other and also different from the BVND/UVND pair. Their computing speeds are relatively close to one another (though a stricter test suggests a slight difference), but distinct from the speeds of BVND and UVND. This indicates clear trade-offs, with no single algorithm definitively outperforming others on both criteria based solely on statistical significance, necessitating an examination of actual performance values (mean OFVs and times) and problem-specific priorities to select the most suitable algorithm.

**Table 5.** Results of Friedman's Test, Nemenyi's Post-Hoc Test, and Pairwise Wilcoxon Test on the OFV of 5-LFL instances.

**Friedman's Test for BVND, PVND, CVND, and UVND**

Number of complete experiments (blocks/rows): 162
Number of algorithms (groups/columns): 4
Friedman chi-squared statistic: 370.926
P-value: 4.387E-80

Note: The Friedman test was conducted to determine if there were any statistically significant differences in the median Objective Function Values (OFVs) achieved by the four algorithms (BVND, PVND, CVND, and UVND), specifically on the 162 5-LFL problem instances. The resulting high chi-squared statistic of 370.926 and an extremely low p-value (4.387E-80), well below the 0.05 significance level, led to the rejection of the null hypothesis that all algorithms produce similar median OFVs. This indicates significant variations in OFV performance among the algorithms when applied to the 5-LFL problem set, justifying further post-hoc analysis to identify which specific algorithms differ in their OFV results.

**--- Post-Hoc Analysis ---**

**Nemenyi's Post-Hoc Test Results (pairwise p-values):**

| | BVND | PVND | CVND | UVND |
|---|---|---|---|---|
| BVND | | 7.0E-06 | 0 | 3.0E-05 |
| PVND | 7.0E-06 | | 0 | 0 |
| CVND | 0 | 0 | | 0 |
| UVND | 3.0E-05 | 0 | 0 | |

Note: Following the significant Friedman test, Nemenyi's post-hoc test was used to perform all pairwise comparisons of the algorithms' OFV performance on the 5-LFL problems, using an alpha of 0.05. The p-values for all comparisons (e.g., BVND vs. PVND, p=7.0E-06; PVND vs. CVND, p=0.0) were well below this threshold. This signifies that statistically significant differences in the OFVs achieved exist between every pair of algorithms. Consequently, for the 5-LFL problems, each algorithm—BVND, PVND, CVND, and UVND—demonstrates a level of OFV performance that is statistically distinct from all the others.

**Pairwise Wilcoxon Test with Bonferroni Correction (p-values)**

| | BVND | PVND | CVND | UVND |
|---|---|---|---|---|
| BVND | | 8.48E-07 | 5.10E-27 | 1.43E-14 |
| PVND | 8.48E-07 | | 2.44E-27 | 6.03E-20 |
| CVND | 5.10E-27 | 2.44E-27 | | 6.73E-27 |
| UVND | 1.43E-14 | 6.03E-20 | 6.73E-27 | |

Note: The pairwise Wilcoxon signed-rank tests, with a Bonferroni correction applied to control the family-wise error rate at 0.05, were also conducted to compare the OFV performance between each pair of algorithms on the 5-LFL problems. These results robustly confirmed Nemenyi's findings, as all Bonferroni-adjusted p-values (e.g., BVND vs. PVND, p=8.478E-07; CVND vs. UVND, p=6.726E-27) were exceptionally small and significantly below the 0.05 alpha level. This provides strong, statistically significant evidence that OFV differences exist between all pairs of the four algorithms, reinforcing that each algorithm achieves a unique distribution of OFVs when tackling the 5-LFL problem set.

**Table 6.** Results of Friedman's Test, Nemenyi's Post-Hoc Test, and Pairwise Wilcoxon Test on the computing time of 5-LFL instances.

**Friedman's Test for BVND, PVND, CVND, UVND**

Number of complete experiments (blocks/rows): 162
Number of algorithms (groups/columns): 4
Friedman chi-squared statistics: 334.778
P-value: 2.948E-72

Note: The Friedman test was applied to assess whether there were statistically significant differences in the median computing times among the four algorithms (BVND, PVND, CVND, and UVND) when solving the 162 *5*-LFL problem instances. With a high chi-squared statistic of 334.778 and an exceptionally low p-value (2.941E-72), far below the 0.05 significance level, the null hypothesis that all algorithms exhibit similar median computing times was decisively rejected. This result strongly indicates significant variations in computational speed among the algorithms for the *5*-LFL problem set, necessitating post-hoc analysis to pinpoint specific pairwise differences in their computing times.

**--- Post-Hoc Analysis ---**

**Nemenyi's Post-Hoc Test Results (pairwise p-values):**

| | BVND | PVND | CVND | UVND |
|---|---|---|---|---|
| BVND | | 0 | 1.70E-07 | 0 |
| PVND | 0 | | 5.51E-05 | 7.98E-14 |
| CVND | 1.70E-07 | 5.51E-05 | | 0 |
| UVND | 0 | 7.98E-14 | 0 | |

Note: Following the Friedman test's indication of overall differences, Nemenyi's post-hoc test was employed to conduct all pairwise comparisons of the algorithms' computing times, specifically for the 5-LFL problems, using an alpha of 0.05. The p-values for all pairs (e.g., BVND vs. PVND, p < 1.0E-07; PVND vs. CVND, p = 5.51E-05) were well below this significance threshold. This outcome demonstrates that statistically significant differences in computing time exist between every single pair of algorithms. Therefore, for the 5-LFL problems, each algorithm—BVND, PVND, CVND, and UVND—operates at a statistically distinct speed from all the others.

**Pairwise Wilcoxon Test with Bonferroni Correction (p-values)**

| | BVND | PVND | CVND | UVND |
|---|---|---|---|---|
| BVND | | 4.73E-26 | 1.47E-07 | 1.53E-26 |
| PVND | 4.73E-26 | | 9.78E-12 | 2.28E-22 |
| CVND | 1.47E-07 | 9.78E-12 | | 5.95E-24 |
| UVND | 1.53E-26 | 2.28E-22 | 5.95E-24 | |

Note: The pairwise Wilcoxon signed-rank tests, adjusted with a Bonferroni correction to maintain a family-wise error rate of 0.05, were also utilized to compare the computing times between each pair of algorithms for the 5-LFL problem set. These more conservative tests strongly corroborated Nemenyi's findings, with all Bonferroni-adjusted p-values (e.g., BVND vs. CVND, p = 1.47E-07; CVND vs. UVND, p = 5.95E-24) extremely small and significantly below the 0.05 alpha level. This provides robust evidence that the computing time for each of the four algorithms is statistically distinguishable from that of every other algorithm when applied to the 5-LFL problems.

When we analyzed the results for the 5-LFL problem set, a distinct performance hierarchy emerged among the four algorithms tested: BVND, PVND, CVND, and UVND. Regarding solution quality, we found that each algorithm produced an Objective Function Value (OFV) that was statistically different from those of the others. The same was true for how fast they ran. Each algorithm settled into a unique, statistically significant computational time when solving these problems. What this complete separation in performance suggests is a classic trade-off. This forces us to look beyond the statistics alone and examine the actual mean OFVs and run times to make a practical decision tailored to the specific demands of a 5-LFL challenge.

## 7. Managerial Implications

Designing a multi-level supply chain is a high-stakes strategic task for retailers like Walmart or Amazon. The choice of optimization method directly affects the trade-off between computational time (decision-making speed) and solution quality (long-term cost savings). Based on our experimental results on large-scale $k$-LFL instances, we provide the following decision framework for supply chain managers and network analysts:

Table 7. Decision of Algorithm Selection

| Scenario / Constraint | Recommended Algorithm | Rationale |
|---|---|---|
| Time-Critical Planning (e.g., rapid scenario analysis, dynamic tactical adjustments) | PVND or CVND | These variants demonstrated significantly lower CPU times while maintaining competitive solution quality. Suitable when "good enough" solutions are needed quickly. |
| Strategic Network Design (e.g., building new infrastructure, 5-10 year planning) | BVND or UVND | For long-term capital investments, solution quality is paramount. These algorithms achieved the lowest Objective Function Values (OFV), justifying the higher computational cost. |
| Highly Dense Networks (High arc density between echelons) | UVND | The Union VND approach explores the solution space more broadly, which is beneficial when the number of feasible connections is vast. |

Managers should first classify their network size. For networks with more than 5,000 customers (Level 1), exact solvers are infeasible. If the fixed costs of facilities (Plants/Warehouses) are very high relative to transportation costs, the BVND approach is preferred as it exhaustively searches the simplest neighborhoods first to optimize facility selection before fine-tuning assignments. We recommend a two-phase implementation. Use PVND to generate a baseline scenario within minutes. If the projected cost exceeds budget constraints, run UVND overnight to squeeze out remaining efficiencies (typically improving OFV by a statistically significant margin as shown in Table 3).

When we think about a company's supply chain, it is easy to picture a simple line from a single factory to a single store. The reality, especially for large retailers, is much more complicated. Their success often depends on designing a multi-level network, which might involve structuring operations across four or even five distinct layers. For a typical retail business, this process could involve moving products from manufacturing plants to massive central warehouses, then to regional distribution centers, and finally onto the shelves of local stores where customers shop. The real puzzle is determining the optimal location and size for each facility, while balancing transportation costs, managing inventory, and operating the retail stores.

The sheer scale of this challenge becomes clear when you consider a giant like Walmart, which serves around 255 million customers weekly and manages over 100,000 products from thousands of suppliers across more than 10,000 stores. Making that system work is not magic; it is the result of decades of careful decisions regarding where to place their plants, warehouses, and distribution centers. These strategic choices have been key in keeping transport costs low and ensuring products are delivered on time.

This complexity can go even deeper. A five-level system is common in e-commerce, where getting a package to someone's front door is everything. Amazon's network, for instance, appears to follow this model: goods move from suppliers to initial sorting hubs, then to large fulfillment centers, to smaller regional sortation centers, and finally to local stations for last-mile delivery. Each step represents its unique puzzle with its own set of trade-offs.

Owning every piece of the puzzle is not the only path to success. While Amazon invests heavily in its distribution infrastructure, apparel companies like Zara and Benetton have succeeded with a different approach. They maintain agility by working with a wide network of small, independent manufacturers. This enables them to manage the first stage of their supply chain, allowing for rapid adaptation to changing fashion trends. Researchers such as Soshko, Merkuryev, and Chakste (2007) have noted that selecting the correct facility location can yield tangible benefits, including lower transportation costs and improved service.

## 8. Conclusion

This study focused on single-assignment MFL problems, specifically for 4- and 5-level locations. All customers (i.e., retail stores) at Level 1 must be served. Single-assignment from upper-level facilities to lower-level facilities is considered. Furthermore, each retail store also prefers products from a specific set of plants. The flow of a bundle of products from an upper-level facility to a lower-level facility incurs some costs each time it is moved. Selecting each facility also incurs a one-time fixed cost, and the number of facilities selected at each level is limited by an upper bound. We considered large-scale problems and provided four variants of VND meta-heuristics. Extensive computational experiments with heuristics are provided for randomly generated problems, and sensitivity analyses, supported by appropriate statistical methods, are used to validate the effectiveness of the heuristics' results.

While this study demonstrates the relative superiority of specific VND variants for large-scale instances, a limitation of the current work is the absence of benchmarks against exact solvers (e.g., CPLEX or Gurobi) for the largest instances due to computational intractability. Future research should focus on establishing lower bounds for these massive $k$-LFL datasets, perhaps using Lagrangian relaxation or Benders decomposition, to provide an absolute optimality gap assessment. Additionally, integrating the proposed heuristics into a commercial solver framework as a warm-start mechanism could bridge the gap between the speed of heuristics and the exact optimality of an exact solver.
Further research may be considered as follows:

- There are different methods to embed sequences within heuristics. In this study, we employed the *l-Opt* strategy adapted from a traveling salesman type application. However, it would be valuable to compare these approaches to determine which performs best for these problems.
- We also used hierarchical problems; however, retailers (also intermediate facilities) often order directly from upper-level facilities. This situation requires further attention, and we are addressing it for future research. In such cases, it is appropriate to consider a multimodal situation, which we are also considering.
- Each time a product bundle is transferred from one facility to the next, it incurs a cost. This cost is independent of the retail store, as in Ortiz-Astorquiza, Contreras, and Laporte (2018). However, it makes sense to consider such costs when they depend on the retail store, as seen in Ortiz-Astorquiza, Contreras, and Laporte (2019). In such cases, the number of variables increases significantly, and the use of computer storage also increases substantially.
- We did not include capacity constraints for the selected facilities. However, practical problems often require capacity consideration, including opening, closing, and expanding facilities. It makes sense to consider such situations for single-assignment problems; however, for multi-assignment problems, such cases have already been considered (Melo, Nickel, and Saldanha da Gama, 2006).

**Disclosure of interest**: The Authors have no conflict of interest

**Computational Results: Table 1(a,b,c,d,e,f) and Table 2 (a,b,c,d,e,f)**

**Table 1a.** Computational results of four VND processes on 4-LFL instances of 2000 retailers.

| **H=High, M=Medium, L=Low, Lg=Large, Med=Medium, Sm=Small (Dens=Density, Fx= Fixed Costs)** | | | | | | | | |
|---|---|---|---|---|---|---|---|---|
| **R=# Retailers, D=# Dist Centers, W=# Warehouses, P=# Plants** | | | | | | | | |
| **Problem ID** | **Objective Function** | | | | **Time to Best (Seconds)** | | | |
| (R-D-W-P-Density-Fixed Cost Size-#) | BVND | PVND | CVND | UVND | BVND | PVND | CVND | UVND |
| 2000-150-50-30-Hdend-LgFx-1 | 82516 | 86936 | 86799 | 82526 | 97.03 | 91.99 | 92.08 | 100.99 |
| 2000-150-50-30-Hdend-LgFx-2 | 93339 | 93795 | 93672 | 93283 | 348.10 | 348.56 | 346.97 | 350.82 |
| 2000-150-50-30-Hdend-LgFx-3 | 104994 | 114192 | 113427 | 105810 | 168.18 | 150.51 | 147.30 | 166.85 |
| 2000-150-50-30-Hdend-MedFx-1 | 87044 | 91414 | 90685 | 87105 | 47.25 | 38.94 | 31.96 | 42.92 |
| 2000-150-50-30-Hdend-MedFx-2 | 77593 | 78426 | 77962 | 77798 | 44.65 | 46.59 | 38.42 | 48.01 |
| 2000-150-50-30-Hdend-MedFx-3 | 84139 | 84258 | 84285 | 84139 | 335.55 | 334.96 | 328.18 | 344.14 |
| 2000-150-50-30-Hdend-SmFx-1 | 95006 | 95295 | 95236 | 94980 | 236.52 | 235.87 | 221.99 | 254.58 |
| 2000-150-50-30-Hdend-SmFx-2 | 88092 | 90839 | 90091 | 88084 | 225.18 | 172.04 | 181.31 | 217.52 |
| 2000-150-50-30-Hdend-SmFx-3 | 79293 | 81267 | 79732 | 79287 | 232.29 | 221.69 | 216.82 | 221.03 |
| 2000-150-50-30-Ldend-LgFx-1 | 48260 | 49314 | 49400 | 48260 | 114.10 | 113.42 | 113.44 | 114.19 |
| 2000-150-50-30-Ldend-LgFx-2 | 53515 | 53675 | 53516 | 53677 | 17.08 | 16.73 | 16.53 | 16.81 |
| 2000-150-50-30-Ldend-LgFx-3 | 49092 | 49289 | 49180 | 49092 | 6.59 | 6.30 | 6.33 | 6.93 |
| 2000-150-50-30-Ldend-MedFx-1 | 43313 | 44401 | 44109 | 43295 | 8.46 | 7.57 | 7.54 | 8.28 |
| 2000-150-50-30-Ldend-MedFx-2 | 44200 | 46591 | 46503 | 44160 | 109.07 | 108.09 | 108.01 | 108.94 |
| 2000-150-50-30-Ldend-MedFx-3 | 42525 | 42774 | 42773 | 42518 | 40.83 | 40.38 | 40.15 | 40.83 |
| 2000-150-50-30-Ldend-SmFx-1 | 37037 | 37132 | 37179 | 37031 | 86.99 | 86.87 | 86.80 | 86.90 |
| 2000-150-50-30-Ldend-SmFx-2 | 41712 | 42716 | 42203 | 41712 | 142.56 | 141.24 | 141.83 | 141.75 |
| 2000-150-50-30-Ldend-SmFx-3 | 43197 | 43268 | 43190 | 43188 | 154.14 | 153.66 | 153.80 | 153.78 |
| 2000-150-50-30-Mdend-LgFx-1 | 87964 | 88780 | 88322 | 87951 | 264.56 | 261.60 | 260.82 | 266.42 |
| 2000-150-50-30-Mdend-LgFx-2 | 86041 | 87853 | 86522 | 86557 | 265.45 | 264.17 | 263.81 | 267.34 |
| 2000-150-50-30-Mdend-LgFx-3 | 96397 | 103287 | 100513 | 96410 | 87.37 | 79.96 | 79.25 | 84.27 |
| 2000-150-50-30-Mdend-MedFx-1 | 78867 | 80260 | 78863 | 78159 | 88.88 | 83.26 | 88.22 | 87.10 |
| 2000-150-50-30-Mdend-MedFx-2 | 76912 | 84061 | 83021 | 76928 | 232.99 | 227.44 | 227.16 | 229.98 |
| 2000-150-50-30-Mdend-MedFx-3 | 69865 | 76446 | 76090 | 69916 | 162.71 | 158.24 | 158.42 | 162.97 |
| 2000-150-50-30-Mdend-SmFx-1 | 63630 | 63711 | 63626 | 63646 | 155.42 | 151.46 | 152.40 | 154.65 |
| 2000-150-50-30-Mdend-SmFx-2 | 83593 | 84916 | 83572 | 83605 | 50.19 | 33.47 | 43.52 | 40.19 |
| 2000-150-50-30-Mdend-SmFx-3 | 71192 | 76370 | 75316 | 71162 | 268.24 | 265.92 | 263.91 | 266.17 |

**Table 1b.** Computational results of four VND processes on 4-LFL instances of 4000 retailers.

| **H=High, M=Medium, L=Low, Lg=Large, Med=Medium, Sm=Small (Dens=Density, Fx= Fixed Costs)** | | | | | | | | |
|---|---|---|---|---|---|---|---|---|
| **R=# Retailers, D=# Dist Centers, W=# Warehouses, P=# Plants** | | | | | | | | |
| **Problem ID** | **Objective Function** | | | | **Time to Best (Seconds)** | | | |
| (R-D-W-P-Density-Fixed Cost Size-#) | BVND | PVND | CVND | UVND | BVND | PVND | CVND | UVND |
| 4000-150-50-30-Hdend-LgFx-1 | 163782 | 167375 | 166115 | 163851 | 1072.31 | 1061.14 | 1062.49 | 1065.11 |
| 4000-150-50-30-Hdend-LgFx-2 | 144186 | 147013 | 146360 | 144151 | 226.74 | 225.45 | 204.67 | 225.98 |
| 4000-150-50-30-Hdend-LgFx-3 | 142331 | 145623 | 144457 | 142388 | 320.24 | 321.90 | 326.89 | 321.63 |
| 4000-150-50-30-Hdend-MedFx-1 | 179338 | 180803 | 179841 | 179462 | 1029.07 | 997.08 | 999.71 | 1003.84 |
| 4000-150-50-30-Hdend-MedFx-2 | 148549 | 148879 | 148540 | 148544 | 268.42 | 252.38 | 251.16 | 254.26 |
| 4000-150-50-30-Hdend-MedFx-3 | 168375 | 176749 | 173930 | 168375 | 1166.66 | 987.03 | 1054.70 | 1034.42 |
| 4000-150-50-30-Hdend-SmFx-1 | 129510 | 129906 | 129520 | 129487 | 569.84 | 542.02 | 546.21 | 566.37 |
| 4000-150-50-30-Hdend-SmFx-2 | 170884 | 174552 | 173076 | 170863 | 317.76 | 204.84 | 235.93 | 212.54 |
| 4000-150-50-30-Hdend-SmFx-3 | 135090 | 136289 | 135533 | 135097 | 222.43 | 176.78 | 183.69 | 187.66 |
| 4000-150-50-30-Ldend-LgFx-1 | 84948 | 88032 | 87999 | 84945 | 604.26 | 601.90 | 601.64 | 604.20 |
| 4000-150-50-30-Ldend-LgFx-2 | 81317 | 82363 | 82189 | 81342 | 139.10 | 137.22 | 136.63 | 138.73 |
| 4000-150-50-30-Ldend-LgFx-3 | 82540 | 82795 | 82644 | 82540 | 694.10 | 692.76 | 692.87 | 694.06 |
| 4000-150-50-30-Ldend-MedFx-1 | 76294 | 77101 | 76816 | 76181 | 472.91 | 472.02 | 470.85 | 472.08 |
| 4000-150-50-30-Ldend-MedFx-2 | 73613 | 74692 | 74580 | 73604 | 134.52 | 134.65 | 133.58 | 134.54 |
| 4000-150-50-30-Ldend-MedFx-3 | 81556 | 82755 | 82031 | 81505 | 752.50 | 750.80 | 751.65 | 751.20 |
| 4000-150-50-30-Ldend-SmFx-1 | 80637 | 81740 | 81145 | 80689 | 710.93 | 708.15 | 708.62 | 708.90 |
| 4000-150-50-30-Ldend-SmFx-2 | 73529 | 74665 | 74167 | 73521 | 522.74 | 519.80 | 520.05 | 521.80 |
| 4000-150-50-30-Ldend-SmFx-3 | 77185 | 79849 | 79605 | 77182 | 260.39 | 258.22 | 258.28 | 259.58 |
| 4000-150-50-30-Mdend-LgFx-1 | 148423 | 159622 | 158510 | 148396 | 478.79 | 458.22 | 460.49 | 464.91 |
| 4000-150-50-30-Mdend-LgFx-2 | 129456 | 135679 | 135126 | 127102 | 238.96 | 221.87 | 223.87 | 228.46 |
| 4000-150-50-30-Mdend-LgFx-3 | 152378 | 163307 | 160447 | 153006 | 1214.70 | 1181.89 | 1187.64 | 1193.33 |
| 4000-150-50-30-Mdend-MedFx-1 | 136152 | 138858 | 137970 | 136156 | 669.24 | 657.77 | 656.69 | 660.70 |
| 4000-150-50-30-Mdend-MedFx-2 | 160547 | 167382 | 165509 | 159933 | 710.64 | 704.40 | 694.18 | 709.02 |
| 4000-150-50-30-Mdend-MedFx-3 | 128098 | 145466 | 143543 | 128109 | 260.81 | 232.19 | 232.36 | 243.18 |
| 4000-150-50-30-Mdend-SmFx-1 | 123517 | 138277 | 136997 | 123519 | 309.29 | 271.82 | 280.25 | 289.08 |
| 4000-150-50-30-Mdend-SmFx-2 | 121252 | 126117 | 124998 | 121237 | 436.82 | 410.81 | 402.68 | 423.95 |
| 4000-150-50-30-Mdend-SmFx-3 | 118064 | 122915 | 121221 | 118063 | 619.87 | 598.78 | 604.78 | 606.84 |

**Table 1c.** Computational results of four VND processes on 4-LFL instances of 5000 retailers.

| **H=High, M=Medium, L=Low, Lg=Large, Med=Medium, Sm=Small (Dens=Density, Fx= Fixed Costs)** | | | | | | | | |
|---|---|---|---|---|---|---|---|---|
| **R=# Retailers, D=# Dist Centers, W=# Warehouses, P=# Plants** | | | | | | | | |
| **Problem ID** | **Objective Function** | | | | **Time to Best (Seconds)** | | | |
| (R-D-W-P-Density-Fixed Cost Size-#) | BVND | PVND | CVND | UVND | BVND | PVND | CVND | UVND |
| 5000-150-50-30-Hdend-LgFx-1 | 165209 | 165790 | 165369 | 165228 | 311.40 | 292.38 | 293.05 | 302.70 |
| 5000-150-50-30-Hdend-LgFx-2 | 212370 | 214486 | 212124 | 211893 | 865.58 | 771.14 | 805.62 | 781.36 |
| 5000-150-50-30-Hdend-LgFx-3 | 147921 | 169689 | 169670 | 147936 | 555.54 | 546.31 | 550.74 | 561.90 |
| 5000-150-50-30-Hdend-MedFx-1 | 201505 | 202249 | 201523 | 201519 | 923.67 | 834.44 | 848.55 | 854.05 |
| 5000-150-50-30-Hdend-MedFx-2 | 168430 | 170054 | 168760 | 168430 | 1795.55 | 1729.24 | 1741.30 | 1758.29 |
| 5000-150-50-30-Hdend-MedFx-3 | 197495 | 198152 | 197498 | 197512 | 1118.70 | 1091.12 | 1096.34 | 1097.77 |
| 5000-150-50-30-Hdend-SmFx-1 | 135116 | 157476 | 157285 | 135116 | 1281.75 | 1237.76 | 1229.08 | 1297.30 |
| 5000-150-50-30-Hdend-SmFx-2 | 205700 | 209432 | 208485 | 205700 | 1075.55 | 995.81 | 985.64 | 1011.30 |
| 5000-150-50-30-Hdend-SmFx-3 | 203097 | 222019 | 218951 | 203101 | 781.38 | 556.89 | 605.73 | 650.90 |
| 5000-150-50-30-Ldend-LgFx-1 | 103976 | 103981 | 103976 | 103976 | 541.60 | 541.73 | 541.59 | 541.53 |
| 5000-150-50-30-Ldend-LgFx-2 | 106521 | 107654 | 106970 | 106481 | 719.53 | 715.77 | 714.87 | 718.08 |
| 5000-150-50-30-Ldend-LgFx-3 | 99511 | 100092 | 99548 | 99451 | 508.95 | 508.38 | 507.21 | 508.10 |
| 5000-150-50-30-Ldend-MedFx-1 | 95877 | 97206 | 96822 | 95885 | 1012.71 | 1011.60 | 1011.04 | 1011.61 |
| 5000-150-50-30-Ldend-MedFx-2 | 88235 | 89074 | 88428 | 88244 | 801.17 | 798.02 | 797.06 | 798.93 |
| 5000-150-50-30-Ldend-MedFx-3 | 96004 | 97554 | 97350 | 96059 | 326.70 | 324.71 | 324.06 | 326.99 |
| 5000-150-50-30-Ldend-SmFx-1 | 90998 | 91652 | 91121 | 91072 | 775.48 | 774.02 | 772.99 | 774.25 |
| 5000-150-50-30-Ldend-SmFx-2 | 92268 | 92616 | 92387 | 92133 | 1054.15 | 1052.52 | 1051.58 | 1053.32 |
| 5000-150-50-30-Ldend-SmFx-3 | 86461 | 87546 | 87086 | 86490 | 553.43 | 551.13 | 550.27 | 552.15 |
| 5000-150-50-30-Mdend-LgFx-1 | 167206 | 183965 | 182350 | 167208 | 1250.08 | 1218.28 | 1227.12 | 1247.28 |
| 5000-150-50-30-Mdend-LgFx-2 | 183179 | 203938 | 200599 | 181107 | 424.98 | 382.73 | 385.31 | 414.35 |
| 5000-150-50-30-Mdend-LgFx-3 | 172385 | 180852 | 179326 | 172432 | 1062.48 | 1028.18 | 1034.57 | 1039.06 |
| 5000-150-50-30-Mdend-MedFx-1 | 177003 | 195466 | 191273 | 178034 | 1119.24 | 1062.35 | 1069.11 | 1072.15 |
| 5000-150-50-30-Mdend-MedFx-2 | 169550 | 180337 | 178808 | 169566 | 1318.56 | 1250.66 | 1272.73 | 1269.32 |
| 5000-150-50-30-Mdend-MedFx-3 | 193137 | 201951 | 200787 | 192227 | 414.68 | 325.98 | 323.76 | 358.58 |
| 5000-150-50-30-Mdend-SmFx-1 | 167334 | 173544 | 172463 | 167497 | 815.40 | 770.28 | 773.20 | 787.58 |
| 5000-150-50-30-Mdend-SmFx-2 | 155656 | 176255 | 174320 | 155625 | 1566.87 | 1485.99 | 1492.47 | 1533.44 |
| 5000-150-50-30-Mdend-SmFx-3 | 149553 | 161147 | 159245 | 149808 | 364.94 | 307.91 | 303.49 | 324.78 |

**Table 1d.** Computational results of four VND processes on 4-LFL instances of 8000 retailers.

| **H=High, M=Medium, L=Low, Lg=Large, Med=Medium, Sm=Small (Dens=Density, Fx= Fixed Costs)** | | | | | | | | |
|---|---|---|---|---|---|---|---|---|
| **R=# Retailers, D=# Dist Centers, W=# Warehouses, P=# Plants** | | | | | | | | |
| **Problem ID** | **Objective Function** | | | | **Time to Best (Seconds)** | | | |
| (R-D-W-P-Density-Fixed Cost Size-#) | BVND | PVND | CVND | UVND | BVND | PVND | CVND | UVND |
| 8000-150-50-30-Hdend-LgFx-1 | 286710 | 287974 | 286970 | 286712 | 6151.77 | 6090.37 | 6092.19 | 6084.82 |
| 8000-150-50-30-Hdend-LgFx-2 | 308767 | 315077 | 311141 | 3091136 | 2943.53 | 2909.41 | 2874.83 | 2965.03 |
| 8000-150-50-30-Hdend-LgFx-3 | 247574 | 255382 | 255250 | 247581 | 1518.00 | 1507.69 | 1504.14 | 1534.80 |
| 8000-150-50-30-Hdend-MedFx-1 | 290102 | 295803 | 293067 | 290112 | 740.91 | 541.49 | 559.88 | 582.32 |
| 8000-150-50-30-Hdend-MedFx-2 | 290469 | 299441 | 297337 | 290606 | 1200.70 | 1020.87 | 1123.96 | 1003.45 |
| 8000-150-50-30-Hdend-MedFx-3 | 328326 | 329411 | 328329 | 328357 | 962.72 | 567.95 | 758.02 | 579.47 |
| 8000-150-50-30-Hdend-SmFx-1 | 272208 | 295780 | 293111 | 272221 | 1000.10 | 404.84 | 493.88 | 569.34 |
| 8000-150-50-30-Hdend-SmFx-2 | 246276 | 246590 | 246354 | 246276 | 1114.02 | 1061.67 | 1070.72 | 1096.07 |
| 8000-150-50-30-Hdend-SmFx-3 | 381493 | 385195 | 382142 | 381510 | 1478.04 | 1262.27 | 1288.34 | 1286.01 |
| 8000-150-50-30-Ldend-LgFx-1 | 158772 | 160923 | 160715 | 158830 | 857.10 | 846.35 | 845.96 | 855.38 |
| 8000-150-50-30-Ldend-LgFx-2 | 143288 | 143440 | 143445 | 143288 | 1147.32 | 1146.80 | 1146.13 | 1147.54 |
| 8000-150-50-30-Ldend-LgFx-3 | 158978 | 160708 | 160287 | 158942 | 613.94 | 603.96 | 607.36 | 608.72 |
| 8000-150-50-30-Ldend-MedFx-1 | 143892 | 145651 | 145237 | 144002 | 3381.92 | 3375.60 | 3372.98 | 3380.55 |
| 8000-150-50-30-Ldend-MedFx-2 | 144136 | 147494 | 147337 | 144093 | 724.30 | 715.60 | 713.60 | 722.83 |
| 8000-150-50-30-Ldend-MedFx-3 | 152146 | 152564 | 152217 | 152115 | 244.76 | 240.96 | 238.78 | 242.61 |
| 8000-150-50-30-Ldend-SmFx-1 | 133503 | 135270 | 134117 | 133575 | 1046.20 | 1039.95 | 1040.04 | 1043.17 |
| 8000-150-50-30-Ldend-SmFx-2 | 140496 | 143090 | 142734 | 140408 | 778.90 | 775.29 | 774.28 | 775.85 |
| 8000-150-50-30-Ldend-SmFx-3 | 142819 | 146484 | 145582 | 142738 | 1457.02 | 1442.55 | 1443.12 | 1454.80 |
| 8000-150-50-30-Mdend-LgFx-1 | 222168 | 247649 | 243531 | 223419 | 4175.72 | 4084.95 | 4097.18 | 4134.04 |
| 8000-150-50-30-Mdend-LgFx-2 | 279576 | 296747 | 294166 | 279657 | 6167.65 | 6080.30 | 6078.07 | 6115.31 |
| 8000-150-50-30-Mdend-LgFx-3 | 285159 | 299181 | 297137 | 285138 | 2710.06 | 2613.99 | 2618.81 | 2654.68 |
| 8000-150-50-30-Mdend-MedFx-1 | 284920 | 298920 | 296919 | 284763 | 1118.66 | 970.73 | 996.28 | 1016.65 |
| 8000-150-50-30-Mdend-MedFx-2 | 320007 | 341952 | 335765 | 319652 | 5150.24 | 4871.87 | 4926.04 | 4931.52 |
| 8000-150-50-30-Mdend-MedFx-3 | 256300 | 259117 | 257165 | 257062 | 1038.47 | 963.42 | 970.94 | 975.63 |
| 8000-150-50-30-Mdend-SmFx-1 | 239190 | 243659 | 241833 | 239265 | 502.92 | 401.59 | 419.89 | 423.74 |
| 8000-150-50-30-Mdend-SmFx-2 | 270643 | 285783 | 281573 | 270628 | 2032.92 | 1867.06 | 1898.94 | 1898.20 |
| 8000-150-50-30-Mdend-SmFx-3 | 255149 | 264036 | 260899 | 255142 | 4774.45 | 4634.19 | 4687.79 | 4657.68 |

**Table 1e.** Computational results of four VND processes on 4-LFL instances of 9000 retailers.

| **H=High, M=Medium, L=Low, Lg=Large, Med=Medium, Sm=Small (Dens=Density, Fx= Fixed Costs)** | | | | | | | | |
|---|---|---|---|---|---|---|---|---|
| **R=# Retailers, D=# Dist Centers, W=# Warehouses, P=# Plants** | | | | | | | | |
| **Problem ID** | **Objective Function** | | | | **Time to Best (Seconds)** | | | |
| (R-D-W-P-Density-Fixed Cost Size-#) | BVND | PVND | CVND | UVND | BVND | PVND | CVND | UVND |
| 9000-150-50-30-Hdend-LgFx-1 | 390442 | 411536 | 406318 | 390132 | 659.46 | 314.92 | 374.93 | 383.77 |
| 9000-150-50-30-Hdend-LgFx-2 | 318904 | 340494 | 339534 | 319112 | 473.23 | 328.34 | 323.12 | 485.47 |
| 9000-150-50-30-Hdend-LgFx-3 | 367271 | 390719 | 387381 | 385050 | 1773.64 | 1576.77 | 1597.60 | 1565.60 |
| 9000-150-50-30-Hdend-MedFx-1 | 349182 | 372814 | 368714 | 349224 | 6672.44 | 6228.07 | 6400.78 | 6308.50 |
| 9000-150-50-30-Hdend-MedFx-2 | 347940 | 356658 | 347862 | 346429 | 3544.22 | 3033.01 | 3297.32 | 3099.96 |
| 9000-150-50-30-Hdend-MedFx-3 | 364524 | 378008 | 372895 | 364516 | 2949.52 | 2585.90 | 2613.23 | 2728.91 |
| 9000-150-50-30-Hdend-SmFx-1 | 246377 | 273784 | 273392 | 246377 | 1034.88 | 860.35 | 853.36 | 991.70 |
| 9000-150-50-30-Hdend-SmFx-2 | 275550 | 278316 | 276224 | 275550 | 4276.80 | 4105.87 | 4124.20 | 4165.33 |
| 9000-150-50-30-Hdend-SmFx-3 | 343131 | 363371 | 351337 | 343118 | 2657.10 | 2114.36 | 2242.05 | 2075.53 |
| 9000-150-50-30-Ldend-LgFx-1 | 161461 | 165181 | 164801 | 161312 | 747.61 | 737.44 | 737.95 | 745.46 |
| 9000-150-50-30-Ldend-LgFx-2 | 184900 | 185579 | 185122 | 184850 | 482.22 | 475.51 | 474.14 | 476.39 |
| 9000-150-50-30-Ldend-LgFx-3 | 166494 | 174804 | 174409 | 166464 | 1698.78 | 1682.08 | 1681.02 | 1692.64 |
| 9000-150-50-30-Ldend-MedFx-1 | 175197 | 175960 | 175653 | 174880 | 287.13 | 279.18 | 275.92 | 281.53 |
| 9000-150-50-30-Ldend-MedFx-2 | 151682 | 153011 | 151811 | 151724 | 2322.94 | 2318.39 | 2315.97 | 2317.93 |
| 9000-150-50-30-Ldend-MedFx-3 | 160917 | 164013 | 163042 | 160970 | 3667.90 | 3657.80 | 3657.32 | 3662.78 |
| 9000-150-50-30-Ldend-SmFx-1 | 149633 | 153113 | 152299 | 149669 | 3851.04 | 3838.30 | 3839.81 | 3843.20 |
| 9000-150-50-30-Ldend-SmFx-2 | 148940 | 154249 | 153723 | 148877 | 2571.52 | 2560.97 | 2557.90 | 2575.92 |
| 9000-150-50-30-Ldend-SmFx-3 | 163000 | 166762 | 166024 | 162976 | 1018.69 | 1003.78 | 1003.06 | 1013.05 |
| 9000-150-50-30-Mdend-LgFx-1 | 302000 | 347963 | 344931 | 301974 | 4782.90 | 4588.68 | 4637.00 | 4680.38 |
| 9000-150-50-30-Mdend-LgFx-2 | 310059 | 327107 | 324182 | 310220 | 708.31 | 597.94 | 593.89 | 634.69 |
| 9000-150-50-30-Mdend-LgFx-3 | 334976 | 351967 | 348056 | 335199 | 2614.49 | 2511.94 | 2509.76 | 2531.13 |
| 9000-150-50-30-Mdend-MedFx-1 | 278898 | 303181 | 301665 | 279041 | 3753.53 | 3590.35 | 3608.75 | 3667.71 |
| 9000-150-50-30-Mdend-MedFx-2 | 257541 | 263062 | 262050 | 257544 | 1474.90 | 1377.54 | 1378.82 | 1391.96 |
| 9000-150-50-30-Mdend-MedFx-3 | 290489 | 359204 | 356493 | 298982 | 5992.55 | 5768.64 | 5767.68 | 5880.00 |
| 9000-150-50-30-Mdend-SmFx-1 | 265834 | 285834 | 281615 | 266086 | 1783.06 | 1589.02 | 1646.48 | 1637.20 |
| 9000-150-50-30-Mdend-SmFx-2 | 282574 | 309572 | 305118 | 282607 | 1253.98 | 1102.05 | 1125.95 | 1157.28 |
| 9000-150-50-30-Mdend-SmFx-3 | 342723 | 361822 | 357501 | 342703 | 4649.82 | 4249.38 | 4384.29 | 4311.06 |

**Table 1f.** Computational results of four VND processes on 4-LFL instances of 10000 retailers.

| **H=High, M=Medium, L=Low, Lg=Large, Med=Medium, Sm=Small (Dens=Density, Fx= Fixed Costs)** | | | | | | | | |
|---|---|---|---|---|---|---|---|---|
| **R=# Retailers, D=# Dist Centers, W=# Warehouses, P=# Plants** | | | | | | | | |
| **Problem ID** | **Objective Function** | | | | **Time to Best (Seconds)** | | | |
| (R-D-W-P-Density-Fixed Cost Size-#) | BVND | PVND | CVND | UVND | BVND | PVND | CVND | UVND |
| 10000-150-50-30-Hdend-LgFx-1 | 321566 | 323926 | 322505 | 321815 | 639.67 | 551.23 | 497.39 | 578.27 |
| 10000-150-50-30-Hdend-LgFx-2 | 379082 | 383488 | 383735 | 379594 | 1690.76 | 1320.62 | 1313.88 | 1348.39 |
| 10000-150-50-30-Hdend-LgFx-3 | 320632 | 352198 | 351363 | 319578 | 1406.69 | 1263.92 | 1221.20 | 1381.74 |
| 10000-150-50-30-Hdend-MedFx-1 | 361291 | 370279 | 369371 | 361294 | 2303.12 | 2204.46 | 2190.18 | 2262.41 |
| 10000-150-50-30-Hdend-MedFx-2 | 335514 | 339266 | 337369 | 335348 | 2765.08 | 2551.57 | 2614.85 | 2625.55 |
| 10000-150-50-30-Hdend-MedFx-3 | 321362 | 333568 | 327664 | 321329 | 921.41 | 729.91 | 813.82 | 791.89 |
| 10000-150-50-30-Hdend-SmFx-1 | 332991 | 345289 | 339160 | 333045 | 4372.67 | 4054.38 | 4104.60 | 4066.57 |
| 10000-150-50-30-Hdend-SmFx-2 | 359455 | 370686 | 361198 | 359384 | 5858.61 | 5321.47 | 5482.52 | 5397.01 |
| 10000-150-50-30-Hdend-SmFx-3 | 330040 | 337091 | 334763 | 330040 | 2824.08 | 2591.13 | 2652.53 | 2610.52 |
| 10000-150-50-30-Ldend-LgFx-1 | 208777 | 209883 | 209766 | 208743 | 3899.71 | 3891.78 | 3889.88 | 3896.27 |
| 10000-150-50-30-Ldend-LgFx-2 | 179812 | 180655 | 180744 | 179846 | 546.14 | 538.99 | 534.02 | 542.51 |
| 10000-150-50-30-Ldend-LgFx-3 | 195325 | 198598 | 198030 | 195322 | 3219.56 | 3207.65 | 3206.80 | 3212.88 |
| 10000-150-50-30-Ldend-MedFx-1 | 167132 | 168517 | 168252 | 167118 | 3514.88 | 3502.36 | 3500.17 | 3510.79 |
| 10000-150-50-30-Ldend-MedFx-2 | 191287 | 195498 | 193375 | 191238 | 3860.14 | 3843.09 | 3839.75 | 3850.53 |
| 10000-150-50-30-Ldend-MedFx-3 | 187463 | 193603 | 193471 | 187488 | 910.96 | 895.64 | 894.97 | 908.78 |
| 10000-150-50-30-Ldend-SmFx-1 | 180076 | 184180 | 185628 | 180094 | 3182.96 | 3173.86 | 3159.89 | 3179.92 |
| 10000-150-50-30-Ldend-SmFx-2 | 166020 | 173044 | 172237 | 166057 | 3895.22 | 3864.76 | 3863.71 | 3883.09 |
| 10000-150-50-30-Ldend-SmFx-3 | 171613 | 175725 | 172520 | 171893 | 2312.89 | 2299.34 | 2298.64 | 2299.03 |
| 10000-150-50-30-Mdend-LgFx-1 | 362394 | 372857 | 368330 | 362641 | 4808.75 | 4684.20 | 4721.71 | 4737.26 |
| 10000-150-50-30-Mdend-LgFx-2 | 328630 | 340188 | 336855 | 328355 | 190.70 | 81.60 | 82.57 | 110.69 |
| 10000-150-50-30-Mdend-LgFx-3 | 362394 | 372857 | 368330 | 362641 | 4380.88 | 4248.51 | 4287.62 | 4302.37 |
| 10000-150-50-30-Mdend-MedFx-1 | 369488 | 377681 | 374475 | 369511 | 1060.09 | 782.66 | 798.33 | 808.85 |
| 10000-150-50-30-Mdend-MedFx-2 | 332321 | 334136 | 332321 | 332321 | 8630.42 | 8529.22 | 8551.88 | 8530.94 |
| 10000-150-50-30-Mdend-MedFx-3 | 362587 | 368086 | 364583 | 362516 | 7527.28 | 7279.72 | 7295.86 | 7301.94 |
| 10000-150-50-30-Mdend-SmFx-1 | 308897 | 319490 | 316897 | 308894 | 5385.31 | 5228.19 | 5225.23 | 5262.53 |
| 10000-150-50-30-Mdend-SmFx-2 | 343946 | 345068 | 343957 | 343961 | 6257.47 | 6040.22 | 6019.37 | 6054.05 |
| 10000-150-50-30-Mdend-SmFx-3 | 370793 | 388874 | 381709 | 370787 | 1034.55 | 512.71 | 586.47 | 576.04 |

**Table 2a.** Computational results of four VND processes on 5-LFL instances of 2000 retailers.

| **H=High, M=Medium, L=Low, Lg=Large, Med=Medium, Sm=Small (Dens=Density, Fx= Fixed Costs)** | | | | | | | | |
|---|---|---|---|---|---|---|---|---|
| **R=# Retailers, D=# Dist Centers, W=# Warehouses, P=# Plants, S=# Suppliers** | | | | | | | | |
| **Problem ID** | **Objective Function** | | | | **Time to Best (Seconds)** | | | |
| (R-D-W-P-Density-Fixed Cost Size-#) | BVND | PVND | CVND | UVND | BVND | PVND | CVND | UVND |
| 2000-150-50-50-100-Hdend-LgFx-1 | 127180 | 129392 | 133848 | 127020 | 1532.10 | 1819.50 | 2328.41 | 2697.78 |
| 2000-150-50-50-100-Hdend-LgFx-2 | 131016 | 147501 | 156709 | 130646 | 1048.17 | 1233.28 | 2138.57 | 1891.16 |
| 2000-150-50-50-100-Hdend-LgFx-3 | 129708 | 130546 | 139464 | 127778 | 725.87 | 2495.46 | 5133.46 | 4187.80 |
| 2000-150-50-50-100-Hdend-MedFx-1 | 119217 | 119399 | 125171 | 119172 | 1374.84 | 2702.95 | 2142.58 | 4721.81 |
| 2000-150-50-50-100-Hdend-MedFx-2 | 104083 | 105558 | 114242 | 102157 | 1693.83 | 1874.25 | 1833.02 | 2233.02 |
| 2000-150-50-50-100-Hdend-MedFx-3 | 108347 | 105161 | 123306 | 102029 | 1619.39 | 3080.19 | 2080.40 | 3854.57 |
| 2000-150-50-50-100-Hdend-SmFx-1 | 88575 | 90635 | 110244 | 87920 | 989.47 | 1135.67 | 1978.85 | 2397.26 |
| 2000-150-50-50-100-Hdend-SmFx-2 | 103840 | 106741 | 119004 | 105105 | 1231.82 | 1501.53 | 1884.11 | 2628.80 |
| 2000-150-50-50-100-Hdend-SmFx-3 | 84333 | 82403 | 85694 | 82589 | 864.18 | 1025.16 | 2691.49 | 1676.86 |
| 2000-150-50-50-100-Ldend-LgFx-1 | 123531 | 128472 | 128979 | 111671 | 331.52 | 234.22 | 507.87 | 1138.07 |
| 2000-150-50-50-100-Ldend-LgFx-2 | 127847 | 130257 | 133605 | 127971 | 638.95 | 616.27 | 676.31 | 793.06 |
| 2000-150-50-50-100-Ldend-LgFx-3 | 128693 | 128809 | 131974 | 125444 | 463.33 | 506.85 | 588.02 | 1010.13 |
| 2000-150-50-50-100-Ldend-MedFx-1 | 102544 | 103234 | 115655 | 96479 | 493.41 | 819.55 | 777.84 | 1223.59 |
| 2000-150-50-50-100-Ldend-MedFx-2 | 102246 | 102280 | 119321 | 102241 | 367.36 | 924.08 | 819.71 | 1367.53 |
| 2000-150-50-50-100-Ldend-MedFx-3 | 100115 | 100740 | 112246 | 90659 | 288.88 | 660.87 | 573.91 | 772.63 |
| 2000-150-50-50-100-Ldend-SmFx-1 | 103094 | 107937 | 112919 | 104700 | 557.78 | 936.14 | 665.24 | 1212.30 |
| 2000-150-50-50-100-Ldend-SmFx-2 | 86687 | 90362 | 114823 | 85569 | 923.04 | 1330.40 | 1022.91 | 1797.37 |
| 2000-150-50-50-100-Ldend-SmFx-3 | 87611 | 88629 | 93278 | 82495 | 504.37 | 748.80 | 702.74 | 910.61 |
| 2000-150-50-50-100-Mdend-LgFx-1 | 133450 | 134299 | 142980 | 134235 | 1050.54 | 927.73 | 1296.05 | 2541.49 |
| 2000-150-50-50-100-Mdend-LgFx-2 | 130008 | 140061 | 141766 | 129883 | 765.22 | 827.62 | 1210.80 | 1094.91 |
| 2000-150-50-50-100-Mdend-LgFx-3 | 141132 | 142844 | 148292 | 141220 | 1223.29 | 1487.47 | 1498.24 | 2290.16 |
| 2000-150-50-50-100-Mdend-MedFx-1 | 107333 | 109648 | 116862 | 106968 | 1161.83 | 2052.24 | 1418.20 | 2751.46 |
| 2000-150-50-50-100-Mdend-MedFx-2 | 115833 | 115557 | 117425 | 109747 | 536.47 | 673.40 | 983.72 | 822.49 |
| 2000-150-50-50-100-Mdend-MedFx-3 | 119481 | 108682 | 126895 | 115650 | 741.21 | 1290.15 | 1455.34 | 2400.86 |
| 2000-150-50-50-100-Mdend-SmFx-1 | 91442 | 93979 | 111419 | 90700 | 1143.57 | 1824.61 | 1423.73 | 2659.35 |
| 2000-150-50-50-100-Mdend-SmFx-2 | 89198 | 91327 | 103848 | 86695 | 309.41 | 625.61 | 643.16 | 886.67 |
| 2000-150-50-50-100-Mdend-SmFx-3 | 80569 | 84387 | 101457 | 81435 | 578.40 | 1337.62 | 1148.89 | 1465.98 |

**Table 2b.** Computational results of four VND processes on 5-LFL instances of 3000 retailers.

| H=High, M=Medium, L=Low, Lg=Large, Med=Medium, Sm=Small (Dens=Density, Fx= Fixed Costs) | | | | | | | | |
|---|---|---|---|---|---|---|---|---|
| **R=# Retailers, D=# Dist Centers, W=# Warehouses, P=# Plants, S=# Suppliers** | | | | | | | | |
| **Problem ID** | **Objective Function** | | | | **Time to Best (Seconds)** | | | |
| (R-D-W-P-Density-Fixed Cost Size-#) | BVND | PVND | CVND | UVND | BVND | PVND | CVND | UVND |
| 3000-150-50-50-100-Hdend-LgFx-1 | 177143 | 181279 | 188765 | 177171 | 1730.51 | 2774.29 | 2855.76 | 3317.83 |
| 3000-150-50-50-100-Hdend-LgFx-2 | 185687 | 181059 | 216072 | 185317 | 2557.96 | 3712.22 | 3452.67 | 4547.93 |
| 3000-150-50-50-100-Hdend-LgFx-3 | 226165 | 229090 | 234020 | 209822 | 1575.09 | 1718.28 | 2555.89 | 6453.75 |
| 3000-150-50-50-100-Hdend-MedFx-1 | 155601 | 143283 | 174113 | 151073 | 3102.29 | 5227.15 | 3224.85 | 8061.98 |
| 3000-150-50-50-100-Hdend-MedFx-2 | 162708 | 181873 | 206139 | 170701 | 1801.24 | 4322.82 | 3228.67 | 6592.10 |
| 3000-150-50-50-100-Hdend-MedFx-3 | 165688 | 168593 | 175312 | 165719 | 2594.15 | 3338.00 | 3202.92 | 5285.76 |
| 3000-150-50-50-100-Hdend-SmFx-1 | 148823 | 154034 | 161306 | 141953 | 2054.82 | 2787.82 | 5039.54 | 5244.57 |
| 3000-150-50-50-100-Hdend-SmFx-2 | 125141 | 127389 | 137693 | 125123 | 1748.31 | 2596.58 | 3565.00 | 4622.17 |
| 3000-150-50-50-100-Hdend-SmFx-3 | 119312 | 120421 | 139086 | 118029 | 2051.55 | 2954.15 | 3483.59 | 4265.77 |
| 3000-150-50-50-100-Ldend-LgFx-1 | 169606 | 181739 | 188606 | 169814 | 1155.88 | 1211.17 | 1245.18 | 1586.68 |
| 3000-150-50-50-100-Ldend-LgFx-2 | 183923 | 181877 | 196333 | 172949 | 1596.45 | 1958.97 | 2159.86 | 2237.93 |
| 3000-150-50-50-100-Ldend-LgFx-3 | 165942 | 169482 | 191824 | 165498 | 1603.40 | 2130.76 | 1663.86 | 2826.26 |
| 3000-150-50-50-100-Ldend-MedFx-1 | 143168 | 133492 | 158848 | 129467 | 994.58 | 1947.49 | 975.56 | 2355.73 |
| 3000-150-50-50-100-Ldend-MedFx-2 | 152122 | 150776 | 154927 | 138277 | 608.50 | 1178.65 | 728.49 | 1612.80 |
| 3000-150-50-50-100-Ldend-MedFx-3 | 139154 | 140911 | 154228 | 134399 | 1636.26 | 1975.74 | 1460.27 | 2385.93 |
| 3000-150-50-50-100-Ldend-SmFx-1 | 149072 | 151165 | 154243 | 144555 | 1489.77 | 2823.59 | 1566.15 | 4249.35 |
| 3000-150-50-50-100-Ldend-SmFx-2 | 131108 | 137822 | 151672 | 131870 | 797.18 | 1442.17 | 851.93 | 1994.85 |
| 3000-150-50-50-100-Ldend-SmFx-3 | 123340 | 129105 | 163749 | 118981 | 814.33 | 1595.61 | 983.17 | 1933.13 |
| 3000-150-50-50-100-Mdend-LgFx-1 | 201766 | 202385 | 219455 | 201890 | 1147.99 | 1695.32 | 1743.13 | 2189.42 |
| 3000-150-50-50-100-Mdend-LgFx-2 | 200961 | 202491 | 211892 | 201378 | 2878.33 | 3390.62 | 3622.60 | 4608.55 |
| 3000-150-50-50-100-Mdend-LgFx-3 | 210361 | 209517 | 220199 | 210284 | 1136.45 | 2181.45 | 2158.88 | 2948.50 |
| 3000-150-50-50-100-Mdend-MedFx-1 | 157920 | 153900 | 162390 | 157980 | 1308.13 | 2323.04 | 1405.11 | 2309.58 |
| 3000-150-50-50-100-Mdend-MedFx-2 | 186210 | 190493 | 203165 | 170739 | 1748.24 | 2765.20 | 1384.80 | 3176.68 |
| 3000-150-50-50-100-Mdend-MedFx-3 | 137290 | 139402 | 172374 | 135231 | 1314.72 | 2252.12 | 2039.26 | 3187.55 |
| 3000-150-50-50-100-Mdend-SmFx-1 | 126919 | 128570 | 138260 | 125620 | 2388.45 | 3446.25 | 2865.53 | 4990.77 |
| 3000-150-50-50-100-Mdend-SmFx-2 | 147877 | 149180 | 149946 | 148057 | 734.46 | 1124.48 | 1752.31 | 3194.35 |
| 3000-150-50-50-100-Mdend-SmFx-3 | 148726 | 154574 | 161357 | 147700 | 1921.16 | 3371.77 | 1957.44 | 4269.02 |

**Table 2c.** Computational results of four VND processes on 5-LFL instances of 4000 retailers.

| H=High, M=Medium, L=Low, Lg=Large, Med=Medium, Sm=Small (Dens=Density, Fx= Fixed Costs) | | | | | | | | |
|---|---|---|---|---|---|---|---|---|
| **R=# Retailers, D=# Dist Centers, W=# Warehouses, P=# Plants, S=# Suppliers** | | | | | | | | |
| **Problem ID** | **Objective Function** | | | | **Time to Best (Seconds)** | | | |
| (R-D-W-P-Density-Fixed Cost Size-#) | BVND | PVND | CVND | UVND | BVND | PVND | CVND | UVND |
| 4000-150-50-50-100-Hdend-LgFx-1 | 250322 | 255850 | 293675 | 245746 | 4784.90 | 5895.74 | 5270.13 | 8884.54 |
| 4000-150-50-50-100-Hdend-LgFx-2 | 249491 | 252572 | 264465 | 249410 | 3299.39 | 4922.62 | 5532.07 | 4923.97 |
| 4000-150-50-50-100-Hdend-LgFx-3 | 312969 | 300518 | 293565 | 301348 | 3644.90 | 5780.33 | 6871.76 | 6377.05 |
| 4000-150-50-50-100-Hdend-MedFx-1 | 241640 | 235244 | 252694 | 229504 | 5714.65 | 12607.79 | 7976.43 | 15501.11 |
| 4000-150-50-50-100-Hdend-MedFx-2 | 195071 | 199395 | 228603 | 194101 | 2018.72 | 6839.02 | 4256.44 | 8104.99 |
| 4000-150-50-50-100-Hdend-MedFx-3 | 240578 | 251823 | 259403 | 239200 | 5801.16 | 9854.55 | 7012.48 | 12639.69 |
| 4000-150-50-50-100-Hdend-SmFx-1 | 159181 | 167744 | 194203 | 159469 | 3210.56 | 5498.01 | 3838.97 | 9474.06 |
| 4000-150-50-50-100-Hdend-SmFx-2 | 222856 | 227498 | 234632 | 221096 | 3180.44 | 5331.01 | 5627.42 | 7754.41 |
| 4000-150-50-50-100-Hdend-SmFx-3 | 190958 | 191280 | 193857 | 181928 | 3453.02 | 5982.70 | 5133.84 | 10205.69 |
| 4000-150-50-50-100-Ldend-LgFx-1 | 187355 | 188975 | 212454 | 187193 | 807.58 | 947.11 | 1122.42 | 1474.65 |
| 4000-150-50-50-100-Ldend-LgFx-2 | 200560 | 204187 | 224145 | 200237 | 1452.09 | 2205.28 | 2472.19 | 2479.21 |
| 4000-150-50-50-100-Ldend-LgFx-3 | 209885 | 213684 | 232702 | 210091 | 1827.84 | 2285.93 | 1577.39 | 2560.29 |
| 4000-150-50-50-100-Ldend-MedFx-1 | 192514 | 197223 | 236585 | 183226 | 2156.44 | 2937.75 | 2830.92 | 3437.66 |
| 4000-150-50-50-100-Ldend-MedFx-2 | 220348 | 200020 | 234815 | 190991 | 2063.87 | 2955.56 | 2439.38 | 3674.91 |
| 4000-150-50-50-100-Ldend-MedFx-3 | 206992 | 201915 | 236051 | 192946 | 2895.13 | 4523.59 | 2737.84 | 5101.41 |
| 4000-150-50-50-100-Ldend-SmFx-1 | 164835 | 167148 | 190134 | 162773 | 2762.72 | 4170.86 | 2754.89 | 4547.10 |
| 4000-150-50-50-100-Ldend-SmFx-2 | 164886 | 168507 | 192300 | 164179 | 1015.01 | 2294.03 | 1326.23 | 3166.26 |
| 4000-150-50-50-100-Ldend-SmFx-3 | 174047 | 171087 | 190750 | 171942 | 3171.29 | 4197.75 | 2680.87 | 5961.77 |
| 4000-150-50-50-100-Mdend-LgFx-1 | 207110 | 209461 | 217585 | 207206 | 2822.33 | 2887.98 | 3620.69 | 3908.46 |
| 4000-150-50-50-100-Mdend-LgFx-2 | 259007 | 259314 | 262621 | 257541 | 2311.42 | 2392.37 | 2963.21 | 3560.05 |
| 4000-150-50-50-100-Mdend-LgFx-3 | 205937 | 207658 | 256148 | 202064 | 1974.51 | 3472.08 | 2630.59 | 4477.17 |
| 4000-150-50-50-100-Mdend-MedFx-1 | 227464 | 223909 | 259425 | 212093 | 4851.74 | 10925.84 | 3933.99 | 12634.54 |
| 4000-150-50-50-100-Mdend-MedFx-2 | 199662 | 197979 | 237652 | 189764 | 2623.72 | 3993.71 | 3472.93 | 5343.94 |
| 4000-150-50-50-100-Mdend-MedFx-3 | 237735 | 221130 | 273895 | 230915 | 1570.35 | 4076.80 | 2939.10 | 2555.40 |
| 4000-150-50-50-100-Mdend-SmFx-1 | 158717 | 162377 | 206016 | 159200 | 2670.78 | 4978.39 | 2834.21 | 7641.99 |
| 4000-150-50-50-100-Mdend-SmFx-2 | 150733 | 145308 | 175914 | 141209 | 1957.93 | 2910.38 | 2509.98 | 4787.07 |
| 4000-150-50-50-100-Mdend-SmFx-3 | 178686 | 181847 | 200345 | 179126 | 3147.32 | 4336.45 | 3606.06 | 6628.40 |

**Table 2d.** Computational results of four VND processes on 5-LFL instances of 5000 retailers.

| **H=High, M=Medium, L=Low, Lg=Large, Med=Medium, Sm=Small (Dens=Density, Fx= Fixed Costs)** | | | | | | | | |
|---|---|---|---|---|---|---|---|---|
| **R=# Retailers, D=# Dist Centers, W=# Warehouses, P=# Plants, S=# Suppliers** | | | | | | | | |
| **Problem ID** | **Objective Function** | | | | **Time to Best (Seconds)** | | | |
| (R-D-W-P-Density-Fixed Cost Size-#) | BVND | PVND | CVND | UVND | BVND | PVND | CVND | UVND |
| 5000-150-50-50-100-Hdend-LgFx-1 | 268947 | 276512 | 297434 | 268887 | 3761.39 | 5433.48 | 6809.43 | 5335.40 |
| 5000-150-50-50-100-Hdend-LgFx-2 | 279601 | 287671 | 307109 | 276429 | 5806.98 | 6224.42 | 8866.12 | 10636.43 |
| 5000-150-50-50-100-Hdend-LgFx-3 | 245074 | 245245 | 250646 | 244790 | 2581.77 | 3593.95 | 5116.11 | 4985.79 |
| 5000-150-50-50-100-Hdend-MedFx-1 | 271725 | 273315 | 330750 | 271971 | 9682.94 | 21701.61 | 8731.48 | 21718.32 |
| 5000-150-50-50-100-Hdend-MedFx-2 | 278547 | 281093 | 290523 | 278039 | 8337.45 | 11630.52 | 4723.23 | 15719.49 |
| 5000-150-50-50-100-Hdend-MedFx-3 | 221034 | 202497 | 335217 | 192316 | 4435.51 | 11194.44 | 4382.31 | 12814.57 |
| 5000-150-50-50-100-Hdend-SmFx-1 | 206043 | 219826 | 263292 | 206468 | 3488.80 | 5877.63 | 4492.62 | 7620.25 |
| 5000-150-50-50-100-Hdend-SmFx-2 | 228559 | 241498 | 245041 | 228941 | 4958.44 | 9709.47 | 4698.77 | 11913.70 |
| 5000-150-50-50-100-Hdend-SmFx-3 | 226301 | 242429 | 270806 | 218908 | 5892.36 | 11834.95 | 7856.14 | 16263.07 |
| 5000-150-50-50-100-Ldend-LgFx-1 | 279697 | 289762 | 294586 | 272549 | 2505.48 | 2329.37 | 3252.98 | 3903.00 |
| 5000-150-50-50-100-Ldend-LgFx-2 | 293944 | 282869 | 306636 | 293751 | 4002.84 | 4481.35 | 4636.14 | 5959.44 |
| 5000-150-50-50-100-Ldend-LgFx-3 | 320792 | 327720 | 362256 | 285965 | 3418.14 | 4958.07 | 3255.28 | 6335.93 |
| 5000-150-50-50-100-Ldend-MedFx-1 | 248485 | 251232 | 284234 | 244833 | 2791.12 | 4775.07 | 3204.98 | 5150.11 |
| 5000-150-50-50-100-Ldend-MedFx-2 | 249046 | 232645 | 283852 | 225577 | 2787.85 | 4544.13 | 2813.65 | 5456.14 |
| 5000-150-50-50-100-Ldend-MedFx-3 | 227814 | 212717 | 248618 | 231710 | 3023.70 | 4814.99 | 3081.34 | 4939.87 |
| 5000-150-50-50-100-Ldend-SmFx-1 | 207160 | 212101 | 237825 | 200998 | 4006.85 | 5428.53 | 3834.05 | 6305.15 |
| 5000-150-50-50-100-Ldend-SmFx-2 | 224958 | 226515 | 263278 | 223008 | 2582.18 | 5077.92 | 2149.58 | 5752.48 |
| 5000-150-50-50-100-Ldend-SmFx-3 | 247893 | 247679 | 265141 | 240354 | 4673.47 | 5704.33 | 4289.94 | 5734.73 |
| 5000-150-50-50-100-Mdend-LgFx-1 | 342111 | 346407 | 365131 | 339762 | 3910.67 | 13969.66 | 3288.19 | 14615.31 |
| 5000-150-50-50-100-Mdend-LgFx-2 | 313128 | 316014 | 321843 | 313308 | 4836.85 | 6858.35 | 5865.27 | 8489.50 |
| 5000-150-50-50-100-Mdend-LgFx-3 | 299188 | 296278 | 315600 | 298995 | 4342.93 | 5757.15 | 6738.68 | 5328.56 |
| 5000-150-50-50-100-Mdend-MedFx-1 | 248908 | 240800 | 325411 | 227999 | 5879.65 | 9970.59 | 5492.00 | 10023.12 |
| 5000-150-50-50-100-Mdend-MedFx-2 | 249029 | 258083 | 294296 | 255775 | 6179.88 | 9984.09 | 7496.14 | 12630.46 |
| 5000-150-50-50-100-Mdend-MedFx-3 | 270201 | 273845 | 330514 | 242249 | 2049.76 | 3460.07 | 2739.98 | 5503.08 |
| 5000-150-50-50-100-Mdend-SmFx-1 | 205491 | 209016 | 238292 | 196062 | 4182.37 | 5773.55 | 3916.20 | 7953.54 |
| 5000-150-50-50-100-Mdend-SmFx-2 | 207958 | 215717 | 267263 | 208563 | 3457.46 | 7494.63 | 2792.22 | 9684.29 |
| 5000-150-50-50-100-Mdend-SmFx-3 | 199180 | 212104 | 254877 | 198394 | 2604.59 | 4328.23 | 3750.37 | 5646.27 |

**Table 2e.** Computational results of four VND processes on 5-LFL instances of 6000 retailers.

| **H=High, M=Medium, L=Low, Lg=Large, Med=Medium, Sm=Small (Dens=Density, Fx= Fixed Costs)** | | | | | | | | |
|---|---|---|---|---|---|---|---|---|
| **R=# Retailers, D=# Dist Centers, W=# Warehouses, P=# Plants, S=# Suppliers** | | | | | | | | |
| **Problem ID** | **Objective Function** | | | | **Time to Best (Seconds)** | | | |
| (R-D-W-P-Density-Fixed Cost Size-#) | BVND | PVND | CVND | UVND | BVND | PVND | CVND | UVND |
| 6000-150-50-50-100-Hdend-LgFx-1 | 293008 | 300028 | 319771 | 271939 | 7518.29 | 7752.59 | 8788.95 | 15675.96 |
| 6000-150-50-50-100-Hdend-LgFx-2 | 393936 | 400790 | 433530 | 392230 | 6817.73 | 9325.12 | 7131.65 | 11823.88 |
| 6000-150-50-50-100-Hdend-LgFx-3 | 331602 | 349018 | 354894 | 312272 | 6512.04 | 7133.33 | 8750.93 | 10520.74 |
| 6000-150-50-50-100-Hdend-MedFx-1 | 309084 | 315475 | 339472 | 309267 | 13910.61 | 23165.08 | 10432.06 | 22768.38 |
| 6000-150-50-50-100-Hdend-MedFx-2 | 321685 | 328873 | 370244 | 294956 | 5233.25 | 7595.34 | 5378.63 | 13408.95 |
| 6000-150-50-50-100-Hdend-MedFx-3 | 306137 | 312871 | 344676 | 276489 | 11065.63 | 17483.96 | 10824.81 | 23178.23 |
| 6000-150-50-50-100-Hdend-SmFx-1 | 298367 | 304797 | 320097 | 298511 | 4744.46 | 8680.57 | 9825.08 | 9010.65 |
| 6000-150-50-50-100-Hdend-SmFx-2 | 276283 | 284867 | 309814 | 271278 | 9608.96 | 14602.37 | 11752.54 | 23380.67 |
| 6000-150-50-50-100-Hdend-SmFx-3 | 233352 | 256432 | 294902 | 234922 | 5538.64 | 8062.13 | 12111.97 | 11105.01 |
| 6000-150-50-50-100-Ldend-LgFx-1 | 297336 | 305185 | 320026 | 297697 | 4015.26 | 5390.77 | 2969.21 | 6951.39 |
| 6000-150-50-50-100-Ldend-LgFx-2 | 320940 | 332249 | 389929 | 320309 | 3010.60 | 4941.23 | 3112.37 | 5173.45 |
| 6000-150-50-50-100-Ldend-LgFx-3 | 383004 | 383836 | 419718 | 382987 | 854.67 | 4519.03 | 63.18 | 2642.41 |
| 6000-150-50-50-100-Ldend-MedFx-1 | 271278 | 234922 | 297697 | 320309 | 9608.96 | 5538.64 | 4015.26 | 3010.60 |
| 6000-150-50-50-100-Ldend-MedFx-2 | 248552 | 255533 | 349898 | 245610 | 5953.86 | 9490.38 | 6630.76 | 10870.23 |
| 6000-150-50-50-100-Ldend-MedFx-3 | 286315 | 288274 | 323095 | 258663 | 2285.24 | 3177.12 | 2570.81 | 4689.73 |
| 6000-150-50-50-100-Ldend-SmFx-1 | 239486 | 245162 | 292851 | 230756 | 4734.46 | 8507.50 | 2910.73 | 10239.04 |
| 6000-150-50-50-100-Ldend-SmFx-2 | 271974 | 274030 | 315028 | 267892 | 5862.15 | 9991.21 | 5813.19 | 11675.50 |
| 6000-150-50-50-100-Ldend-SmFx-3 | 266301 | 267334 | 301369 | 253206 | 2144.32 | 4919.35 | 2499.02 | 6916.76 |
| 6000-150-50-50-100-Mdend-LgFx-1 | 360112 | 394691 | 403522 | 359874 | 3986.33 | 7790.15 | 4217.82 | 10627.67 |
| 6000-150-50-50-100-Mdend-LgFx-2 | 377857 | 384473 | 399858 | 344081 | 3123.52 | 2665.16 | 4604.86 | 7212.37 |
| 6000-150-50-50-100-Mdend-LgFx-3 | 321651 | 328403 | 366359 | 322146 | 5392.56 | 6033.74 | 5172.69 | 7447.67 |
| 6000-150-50-50-100-Mdend-MedFx-1 | 324842 | 318890 | 346186 | 289247 | 5945.23 | 10348.75 | 7156.49 | 14427.86 |
| 6000-150-50-50-100-Mdend-MedFx-2 | 318968 | 291904 | 330334 | 282535 | 8925.05 | 13132.62 | 8354.21 | 15052.72 |
| 6000-150-50-50-100-Mdend-MedFx-3 | 307249 | 289231 | 313039 | 299443 | 5874.64 | 15787.05 | 4672.66 | 11380.20 |
| 6000-150-50-50-100-Mdend-SmFx-1 | 289793 | 283423 | 305845 | 279657 | 7945.99 | 15302.42 | 8741.47 | 15444.77 |
| 6000-150-50-50-100-Mdend-SmFx-2 | 216721 | 224010 | 300369 | 216542 | 5189.52 | 10202.04 | 5511.24 | 12777.54 |
| 6000-150-50-50-100-Mdend-SmFx-3 | 303470 | 311863 | 325223 | 303353 | 6194.88 | 13501.47 | 4643.48 | 18602.30 |

**Table 2f.** Computational results of four VND processes on 5-LFL instances of 8000 retailers.

| **H=High, M=Medium, L=Low, Lg=Large, Med=Medium, Sm=Small (Dens=Density, Fx= Fixed Costs)** | | | | | | | | |
|---|---|---|---|---|---|---|---|---|
| **R=# Retailers, D=# Dist Centers, W=# Warehouses, P=# Plants, S=# Suppliers** | | | | | | | | |
| **Problem ID** | **Objective Function** | | | | **Time to Best (Seconds)** | | | |
| (R-D-W-P-Density-Fixed Cost Size-#) | BVND | PVND | CVND | UVND | BVND | PVND | CVND | UVND |
| 8000-150-50-50-100-Hdend-LgFx-1 | 460043 | 512460 | 580678 | 459786 | 5149.86 | 9977.25 | 7275.21 | 12243.60 |
| 8000-150-50-50-100-Hdend-LgFx-2 | 387970 | 400420 | 489917 | 387969 | 8819.33 | 17220.76 | 10960.60 | 12918.66 |
| 8000-150-50-50-100-Hdend-LgFx-3 | 481098 | 483133 | 495038 | 480988 | 6206.11 | 7531.56 | 10709.55 | 10433.51 |
| 8000-150-50-50-100-Hdend-MedFx-1 | 403368 | 412470 | 474116 | 403122 | 16824.32 | 25967.07 | 15196.93 | 28814.66 |
| 8000-150-50-50-100-Hdend-MedFx-2 | 309120 | 325281 | 480938 | 310566 | 16897.69 | 28653.49 | 17992.02 | 31458.59 |
| 8000-150-50-50-100-Hdend-MedFx-3 | 348911 | 328797 | 388588 | 320849 | 6766.82 | 9356.48 | 10818.30 | 12112.75 |
| 8000-150-50-50-100-Hdend-SmFx-1 | 322998 | 328336 | 439910 | 319933 | 7232.40 | 16204.02 | 8777.54 | 18275.64 |
| 8000-150-50-50-100-Hdend-SmFx-2 | 356273 | 361349 | 382524 | 336157 | 9029.71 | 21010.80 | 11851.64 | 24648.14 |
| 8000-150-50-50-100-Hdend-SmFx-3 | 329781 | 333136 | 354270 | 330067 | 9131.13 | 20176.41 | 9936.22 | 24549.03 |
| 8000-150-50-50-100-Ldend-LgFx-1 | 436168 | 439563 | 462480 | 435383 | 6078.60 | 11017.90 | 4954.90 | 15366.45 |
| 8000-150-50-50-100-Ldend-LgFx-2 | 338560 | 351882 | 380284 | 338346 | 3217.08 | 3795.77 | 3264.54 | 5143.09 |
| 8000-150-50-50-100-Ldend-LgFx-3 | 404454 | 408439 | 443446 | 393934 | 6219.02 | 6657.13 | 6497.89 | 9748.62 |
| 8000-150-50-50-100-Ldend-MedFx-1 | 398400 | 399903 | 428327 | 379651 | 9199.84 | 13199.69 | 9416.97 | 13769.74 |
| 8000-150-50-50-100-Ldend-MedFx-2 | 376643 | 380956 | 410772 | 353777 | 8778.36 | 13618.16 | 6957.65 | 15815.32 |
| 8000-150-50-50-100-Ldend-MedFx-3 | 462088 | 468267 | 506080 | 462155 | 7801.79 | 8852.34 | 6081.72 | 14368.71 |
| 8000-150-50-50-100-Ldend-SmFx-1 | 294810 | 304183 | 368753 | 294360 | 8152.99 | 11967.96 | 7991.52 | 12275.05 |
| 8000-150-50-50-100-Ldend-SmFx-2 | 346083 | 342947 | 387215 | 329934 | 7002.97 | 12752.58 | 5589.84 | 12973.57 |
| 8000-150-50-50-100-Ldend-SmFx-3 | 337552 | 349742 | 407553 | 320960 | 5832.42 | 10282.20 | 3876.65 | 12390.54 |
| 8000-150-50-50-100-Mdend-LgFx-1 | 544783 | 559126 | 636860 | 543943 | 9532.67 | 14818.14 | 10711.08 | 18589.10 |
| 8000-150-50-50-100-Mdend-LgFx-2 | 422654 | 429948 | 460233 | 417468 | 6241.39 | 10765.83 | 5411.63 | 12349.77 |
| 8000-150-50-50-100-Mdend-LgFx-3 | 441001 | 442948 | 443705 | 441232 | 5948.38 | 9164.54 | 8584.71 | 11328.28 |
| 8000-150-50-50-100-Mdend-MedFx-1 | 400289 | 375619 | 406227 | 365818 | 2847.80 | 8144.00 | 4882.62 | 9706.31 |
| 8000-150-50-50-100-Mdend-MedFx-2 | 429118 | 426793 | 487408 | 390742 | 8374.05 | 14886.71 | 8762.91 | 25737.91 |
| 8000-150-50-50-100-Mdend-MedFx-3 | 411866 | 414954 | 425643 | 405364 | 9995.88 | 13812.95 | 9548.09 | 16999.94 |
| 8000-150-50-50-100-Mdend-SmFx-1 | 396622 | 403721 | 424452 | 409841 | 11948.44 | 20194.12 | 11669.24 | 27102.25 |
| 8000-150-50-50-100-Mdend-SmFx-2 | 334539 | 339190 | 363120 | 331733 | 11790.59 | 13732.46 | 11633.42 | 19720.01 |
| 8000-150-50-50-100-Mdend-SmFx-3 | 332512 | 340180 | 417354 | 329292 | 16513.92 | 28054.14 | 12889.34 | 28668.03 |

# Appendix A

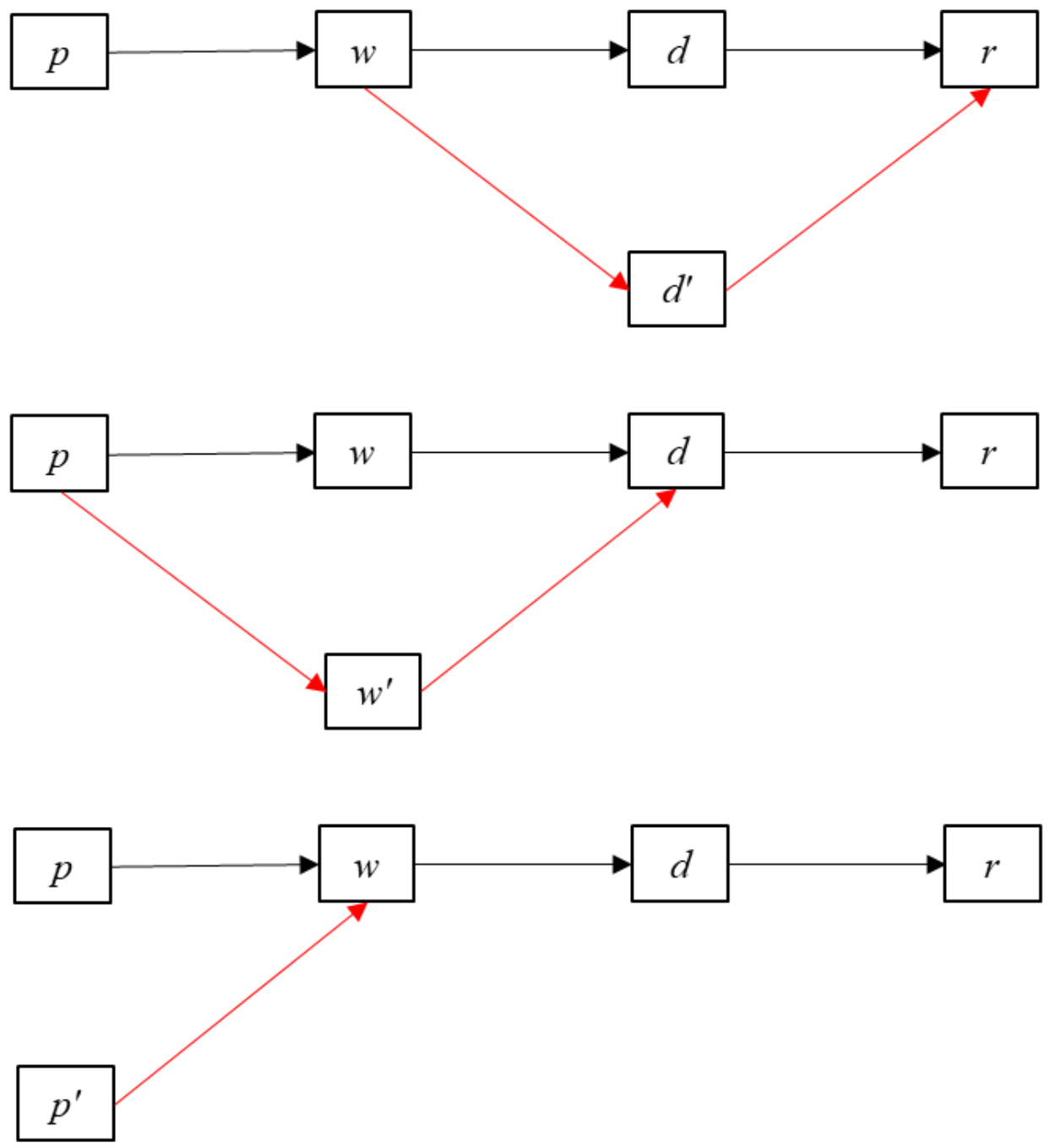

**Figure 2(a)** One-exchnahe, N(1), for 4-LFL

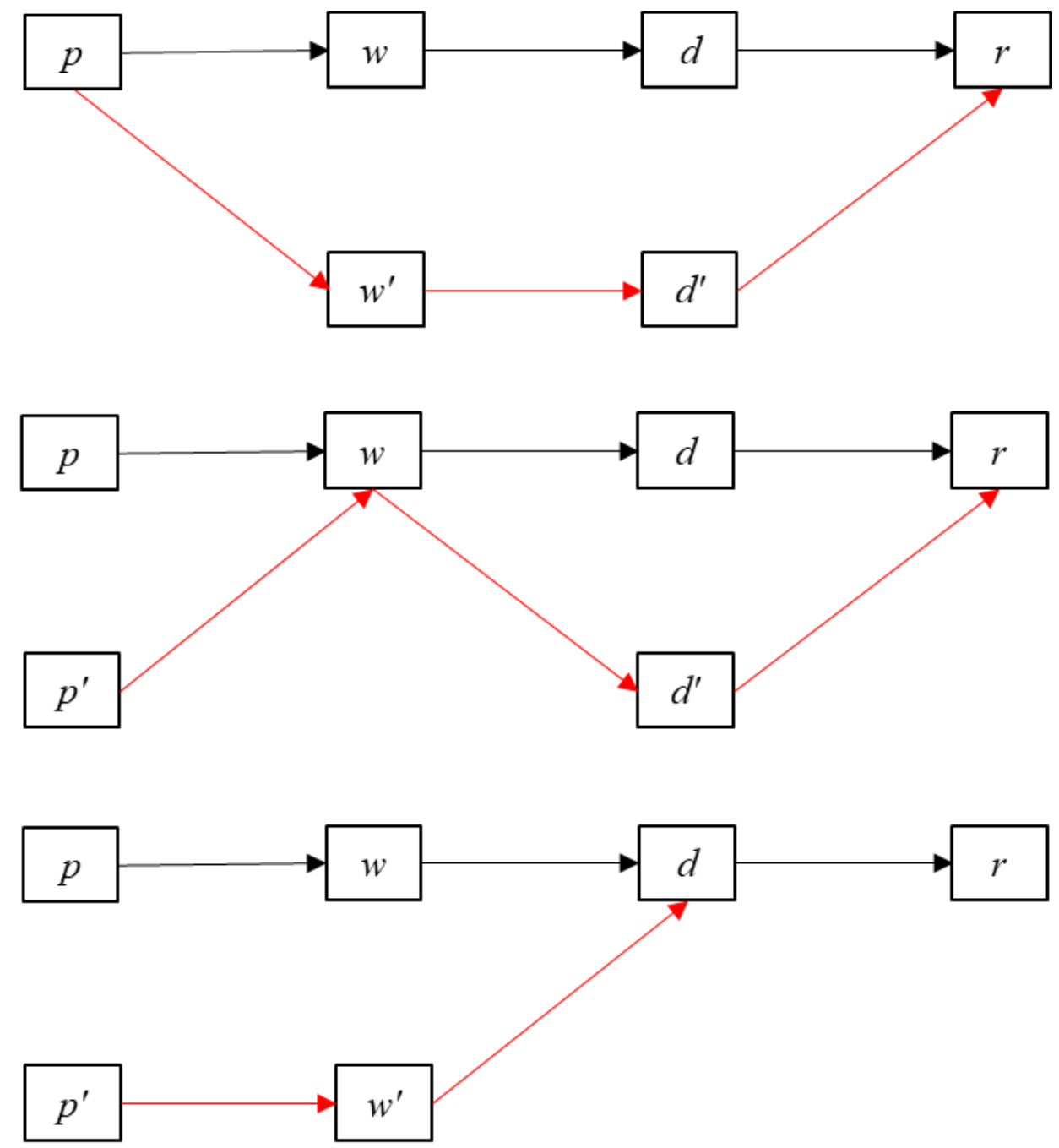

**Figure 2(b)** Two-exchnahe, N(2), for 4-LFL

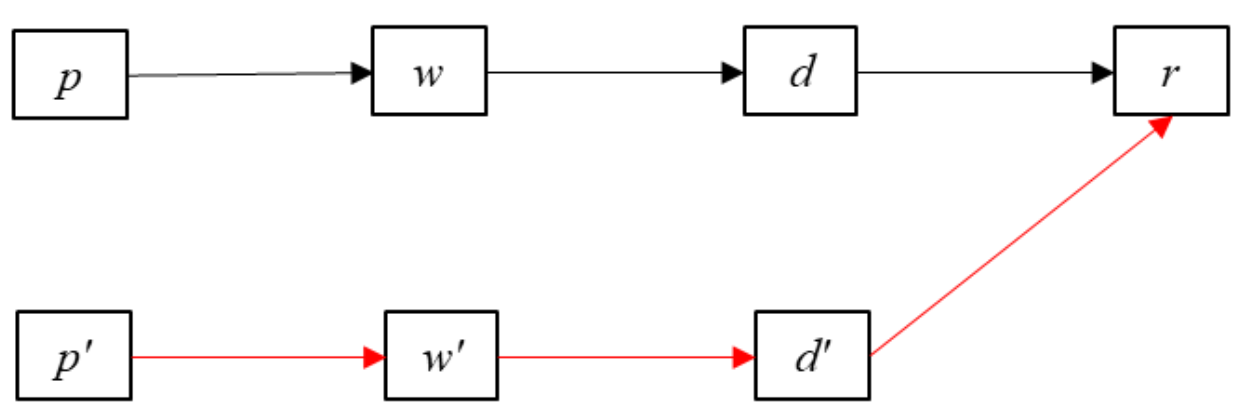

**Figure 2(c)** Three-exchnahe, N(3), for 4-LFL

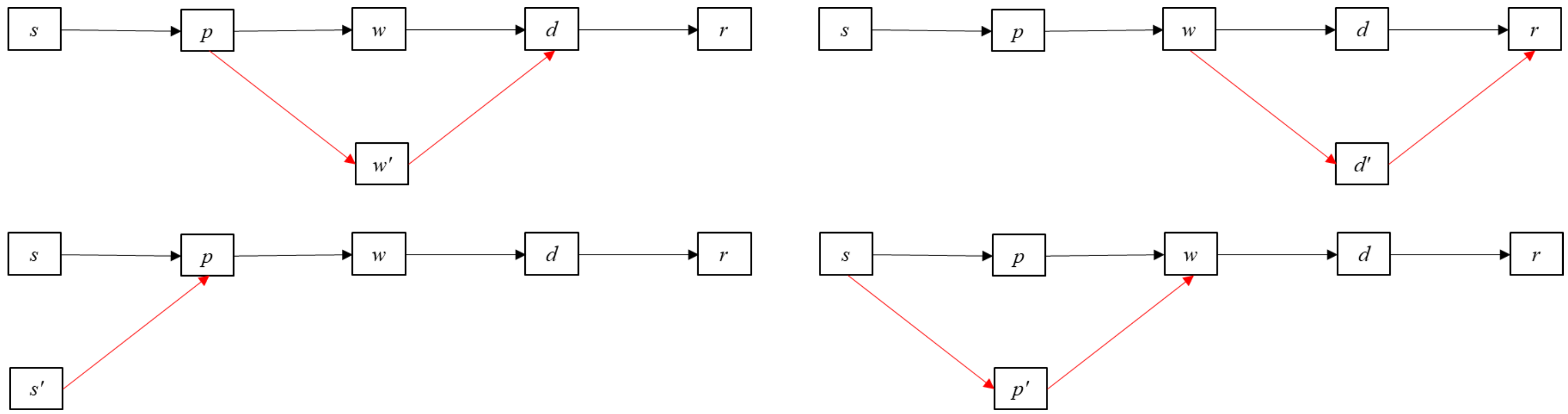

**Figure 3(a)** One-exchnahe, N(1), for 5-LFL

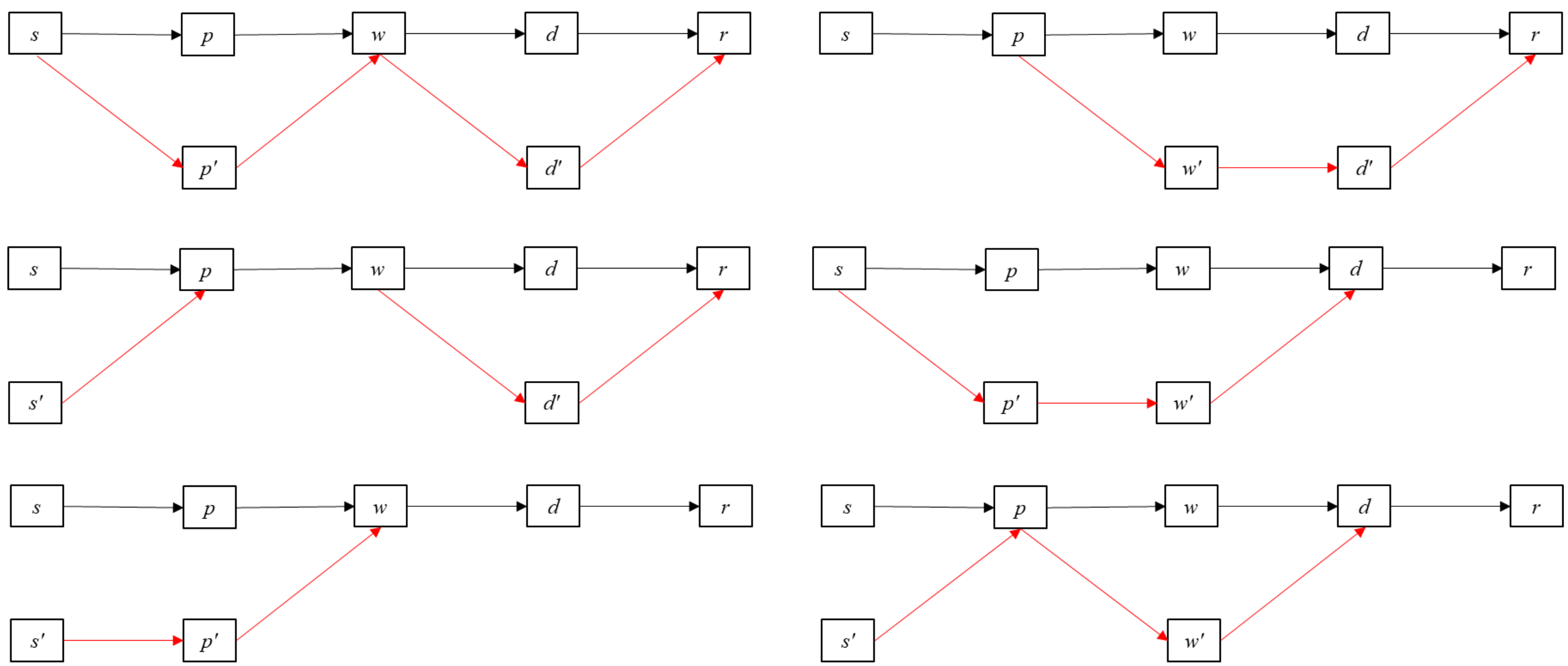

**Figure 3(b)** Two-exchnahe, N(2), for 5-LFL

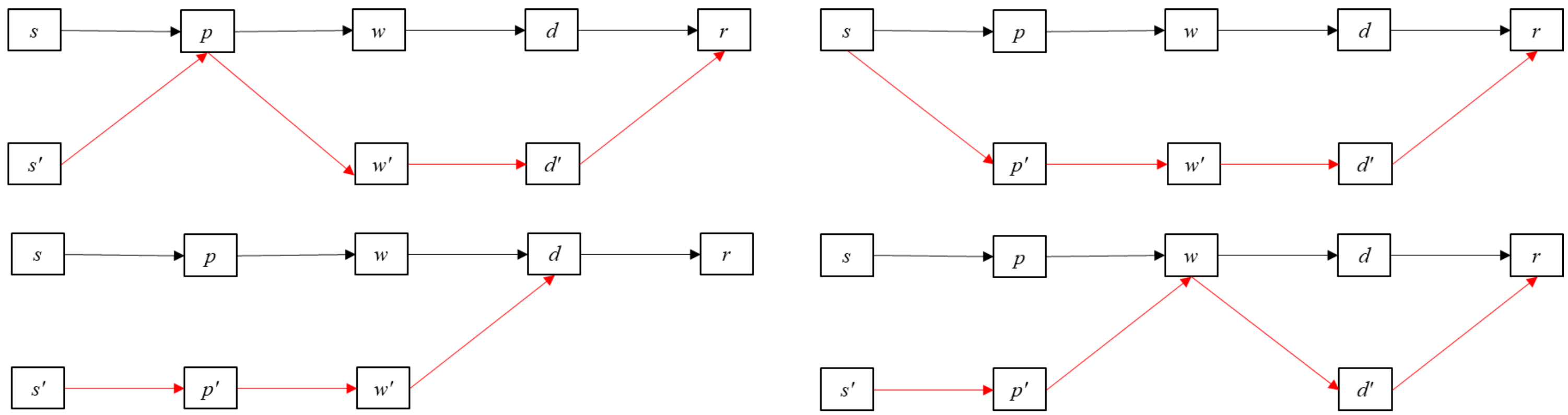

**Figure 3(c**) Three-exchnahe, N(3), for 5-LFL

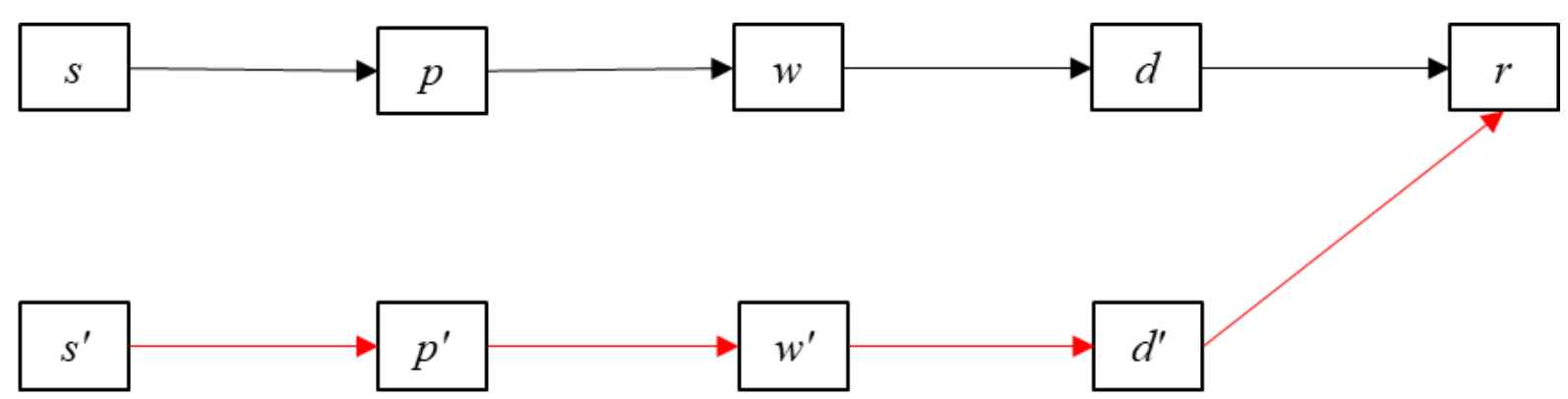

**Figure 3(d**) Four-exchnahe, N(3), for 5-LFL

## Appendix B

Table B1. Parameters used for data generation

| **Parameters:** |
|---|
| **R**, number of retail stores (Level 1), 2000, 3000, 4000, 5000, 6000, 8000, 9000, 10000 |
| For 4-LFL, 2000, 4000, 5000, 8000, 9000, 10000 |
| For 5-LFL, 2000, 3000, 4000, 5000,6000, 8000 |
| **D**, number of distribution centers (Level 2), 150 |
| **W**, number of warehouses (Level 3), 50 |
| **P**, number of plants (Level 4), 30, 50 |
| **S**, number of suppliers (Level 5), 100 |
| **r**, a retail store |
| **d**, a distribution center |
| **w**, a warehouse |
| **p**, a plant |
| **s**, a supplier |
| **P-R, P by R binary matrix**, pr-th element 1 means plant p is eligible to ship products to a retail store r, 0 otherwise. |
| **D-R, D by R matrix**, positive dr-th element means distribution center d is eligible to ship products to retail store r with cost or distance equal to the value of the dr-th element , and 0 means ineligible. |
| **W-D, W by D matrix**, positive wd-th element means warehouse w is eligible to ship products to distribution center d with cost or distance equal to value of wd-th element , and 0 means ineligible. |
| **P-W, P by W matrix**, positive pw-th element means plant p is eligible to ship products to warehouse w with cost or distance equal to the value of pw-th element , and 0 means ineligible. |
| **S-P, S by P matrix**, positive sp-th element means supplier s is eligible to ship products to plant p with cost or distance equal to the value of the sp-th element , and 0 means ineligible. |
| **L_DR (U_DR)**, lower (upper) value for elements of cost matrix D-R, (L_DR, U_DR) = (5, 50) |
| **L_WD (U_WD)**, lower (upper) value for elements of cost matrix W-D, (L_WD, U_WD) = (100, 500) |
| **L_PW (U_PW)**, lower (upper) value for elements of cost matrix P-W, (L_PW, U_PW) = (5, 500) |
| **L_SP (U_SP)**, lower (upper) value for elements of cost matrix S-P, (L_SP, U_SP) = (5, 150) |
| **Ldense**, low density, for 4-FLP 20%, for 5-FLP 40% |
| **Mdense**, medium density, for 4-FLP 40%, for 5-FLP 50% |
| **Hdense**, high density, for 4-FLP 60%, for 5-FLP 60% |
| **L_FD (U_FD)**, lower (upper) value for distribution center fixed costs |
| **L_FW (U_FW)**, lower (upper) value for warehouse fixed costs |
| **L_FP (U_FP)**, lower (upper) value for plant fixed costs |
| **L_FS (U_FS)**, lower (upper) value for supplier fixed costs |
| **SmFx**, small fixed costs. (L_FD, U_FD) = (50, 100), (L_FW, U_FW) = (100, 200), (L_FP,U_FP) = (200, 400), (L_FS, U_FS)= (20, 100) |
| **MedFx**, medium fixed costs, (L_FD, U_FD) = (100, 200), (L_FW, U_FW) = (200, 400), |

(L_FP, U_FP) = (400, 800), (L_FS, U_FS) = (50, 200)
**LgFx**, Large fixed costs, (L_FD, U_FD) = (200, 400), (L_FW, U_FW) = (400, 800), (L_FP, U_FP) = (800, 1600), (L_FS, U_FS) = (200, 400)

---

**UB_D**, upper bound for number of distribution centers to be opened, density*D

**UB_W**, upper bound for number of warehouse centers to be opened, density*W

**UB_P**, upper bound for number of plants to be opened, density*P

**UB_S**, upper bound for number of suppliers to be opened (considered), density*S

---

**Max_Local**, the number of multi-start. For low-density 70, medium-density 50, and high-density 30

---

**Variables:**

**D_Upper**, number of distribution centers opened

**W_Upper**, the number of warehouses opened

**P_Upper**, number of plants opened

**S_Upper**, the number of suppliers opened

**(p, w, d, r)**, schedule of receiving a bundle of products to retailer r, via plant p, warehouse w, and distribution center d in 4-FLP

**(s, p, w, d, r)**, schedule of receiving a bundle of products to retail store r, via supplier s, plant p, warehouse w, and distribution center d in 5-FLP

---

**Objective function:**

Find a schedule of shipments satisfying all retailers, minimizing the total cost of shipment, and opening facilities

## Appendix C

Table C1. Comparison between Gurobi Solver (7200-second time limit) and VND variants (60-second time limit) on 5-LFL instances in terms of objective function value (OFV)

| ID | R-D-W-P-S | Gurobi OFV | OFV of Best VND Variant | GAP % |
|---|---|---|---|---|
| 1 | 100-50-25-25-50 | 8553 | 8726 (BVND) | 2.02% |
| 2 | 100-50-25-25-50 | 8762 | 8762 (BVND) | 0 |
| 3 | 100-50-25-25-50 | 9356 | 9400 (All variants) | 0.47% |
| 4 | 200-50-25-25-50 | 14795 | 15480 (CVND) | 4.63% |
| 5 | 200-50-25-25-50 | 15446 | 15952 (PVND) | 3.28% |
| 6 | 200-50-25-25-50 | 16337 | 17162 (PVND) | 4.79% |
| 7 | 200-50-25-25-50 | 15181 | 16183 (BVND) | 6.60% |